\magnification=1200
\input amstex
\documentstyle{amsppt}

\def\wt{\widetilde}
\def\ssm {\smallsetminus}
\def\SO {\Cal O}
\def\SD {\Cal D}
\def\SE {\Cal E}
\def\SL {\Cal L}
\def\SF {\Cal F}

\def\SU {\Cal U}
\def\SW {\Cal W}

\def\p {\partial}
\topmatter
\title Logarithmic heat projective operators  
\endtitle
\author Xiaotao Sun   \endauthor
\address Institute of Mathematics, Academy of Mathematics and 
Systems Sciences, Chinese Academy of Sciences 
Beijing 100080, China \endaddress
\email xsun$\@$math08.math.ac.cn\endemail
\thanks This research was done when the author was visiting Universit{\"a}t
Essen supported by DFG Forschergruppe "Arithmetik und Geometrie".  \endthanks
\endtopmatter
\document

\heading Introduction\endheading

Let $\Cal C\to S$ be a proper flat family of stable curves, smooth over
$S_0=S\ssm\Delta$. One can associate a flat family 
$f:\Cal {SU}_{\Cal C}(r,d)\to S_0$
of moduli spaces $\Cal {SU}_{\Cal C_s}(r,d)$ of semistable vector bundles of rank $r$ and 
degree $d$ with
fixed determinant over the curves $\Cal C_s$, and also a line bundle $\Theta$
on $\Cal{SU}_{\Cal C}(r.d)$ such that its restriction to each fibre is the line bundle
on $\Cal {SU}_{\Cal C_s}(r,d)$ defined by the theta divisors.
It is well known that for any positive integer $k$ the direct image
$\SE_0:=f_*\Theta^k$ is a vector bundle on $S_0$ and have a flat projective
(Hitchin) connection. A natural question suggested by conformal fields theory
is that: Is the connection regular ? By my knowledge, this problem seems open,
which is part of the motivation of this paper.

On the other hand, let $S$ be a smooth curve, $\Delta=\{s_0\}$ a closed point
and $\Cal C_{s_0}$ be irreducible, smooth but one node $p_0$, we can associate 
a flat family $f:\Cal U_{\Cal C}(r,d)\to S$ of moduli spaces 
$\Cal U_{\Cal C_s}(r,d)$ of semistable torsion free sheaves of rank $r$ and 
degree $d$ over the curves $\Cal C_s$, and also a line bundle $\Theta$
on $\Cal U_{\Cal C}(r.d)$ such that its restriction to each fibre is the 
theta line bundle on $\Cal U_{\Cal C_s}(r,d)$. It was proved in [NR] and [Su]
a factorization of the space $H^0(\Cal U_{\Cal C_s}(r,d),\Theta_{s_0}^k)$,
which is unfortunately not canonical. One may ask for a canonical factorization,
or ask further the questions: if there exists a logarithmic projective 
connection on $\SE=f_*\Theta^k$ such that its restriction to $(\SE_0,S_0)$ is 
the flat projective connection satisfying Hitchin's requirements ? if the answer
is yes, what is the relationship between the residue of the connection and
the factorization one proved in [NR] and [Su] ? In other words, even if we
knew that Hitchin's connection is regular, it is still an interesting problem
to know for what extension of $\SE_0$ such that the connection extends. Thus we
formulate the question as following

\proclaim{Problem} For the family $\SU_{S_0}\to S_0$ of moduli spaces of semistable 
vector bundles and the relative theta line bundle $\Theta_{S_0}$, find the
`correct' degeneration $(\SU_S, \Theta_S)$ of moduli spaces and theta line 
bundles (in other words, the `correct' algeo-geometric analogy of spaces of conformal blocks
on singular curves) such that the direct image of $\Theta^k_S$ is a vector
bundle on $S$ with a flat logarithmic projective connection.\endproclaim    

This paper is a start of the program to study this problem. As the same with
[GJ] and [Hi], we adopt the idea of [We] by considering the existence of heat
operators on a line bundle, which will give the required connection. We first
consider the problem in a general situation (forgeting moduli spaces). Let
$f:X\to S$ be a flat family of (reduced) normal 
crossing varieties of dimension $d$, and $X$, $S$ smooth. 
Assume that $\Delta\subset S$ is a normal crossing divisor and 
$Y:=f^{-1}(\Delta)\subset X$ 
is also a normal crossing divisor such that $f:X\ssm Y\to S\ssm\Delta$ is smooth
and $(X,logY)@>f>>(S,log\Delta)$ is log smooth. Let $\SL$ be a line bundle
on $X$. Then we defined the logarithmic analogies (see Definition 2.2 and 
Definition 2.3) of (projective) heat operators of [GJ] and figured out the
sufficient conditions (see Theorem 2.1) of existence of a projective logarithmic heat 
operator on $\SL$ over $S$, which, similar to [GJ], gave a logarithmic
projective connection on $f_*\SL$. It becomes clear that we need a correct
logarithmic analogy of $sheaf$ $of$ $differential$ $operators$ and have to
work in logarithmic algebraic geometry. We defined the $sheaf$ $of$ $logarithmic$ 
$differential$ $operators$ $on$ $a$ $log$ $scheme$, which may not be a good
definition for general log schemes, but it works well for the log schemes we
concern in this paper (see Proposition 2.1). These materials may be well known
to experts but I am not able to find a reference satisfying my requirements.
Then we checked the conditions in Theorem 2.1 for the rank one case, namely,
we proved that a family of moduli spaces of torsion free sheaves of rank one
over nodal curves satisfies the conditions in Theorem 2.1, thus there
exists a projective logarithmic heat operator on the relative theta line bundle
$\Theta^k$ (see Theorem 3.1).

We developed the necessary technique tools, especially the sheaf of differential
operators in Section 1. Then, in Section 2, we figured out the sufficient 
conditions of existence of a projective logarithmic heat operator (thus the 
conditions of existence of the required projective logarithmic connection) in 
the general situation, and we also gave some descriptions of the conditions.
Finally, in Section 3, we verified the conditions figured out in {\S 2} for
a family of generalized Jacobians (moduli spaces of torsion free sheaves of rank one), 
and thus showed the existence of projective logarithmic heat operator in this case.

\remark{Acknowledgements: This work was done during my stay at FB 6 Mathematik
of Universit{\"a}t Essen. I would like to express my hearty thanks
to Prof. H. Esnault and Prof. E. Viehweg for their hospitality and 
encouragements. I was benefited from the stimulating discussions with them,
which stimulated me to get Lemma 3.3 and Definition 1.3. The discussions
with H. Clemens, I-Hsun Tsai, Kang Zuo, and emails of Z. Ran, are very helpful. 
I thank them very much.}
\endremark

\heading \S1 Logarithmic schemes and logarithmic operators\endheading

In this preliminary section, we recall the so called logarithmic structures
(or log structures for simplicity) on schemes (see [KK]), and define the
sheaves of differential operators on logarithmic schemes. All monoids $M$ are
commutative monoids with unit element and $M^{gp}=\{ab^{-1}\}$ is the associated group.

By a pre-log structure $M$ on a scheme $X$, we mean a sheaf of monoids $M$ on
the {\'e}tale site $X_{et}$ endowed with a homomorphism $\alpha: M\to \SO_X$
with respect to the multiplication on $\SO_X$. A morphism 
$$f:X^{\dagger}:=(X,M)\to Y^{\dagger}:=(Y,N)$$
of schemes with pre-log structures is defined to be a pair $(f,h)$ of a morphism of 
schemes $f:X\to Y$ and a homomorphism $h:f^{-1}(N)\to M$ such that 
$$\CD   
 f^{-1}(N)@>h>>M \\
  @VVV    @VVV     \\
 f^{-1}(\SO_Y) @>>> \SO_X
\endCD$$ 
is commutative. A pre-log structure $(M,\alpha)$ is called a logarithmic structure if 
$$\alpha^{-1}(\SO_X^*)\cong\SO_X^*\quad \text{via} \quad \alpha$$
where $\SO^*_X$ denotes the group of invertible elements of $\SO_X$. A morphism
of schemes with log structures is defined as a morphism of schemes with pre-log
structures. These schemes are called log schemes.

For a pre-log structure $(M,\alpha)$ on $X$, one can define its associated log structure 
$(M^a, \alpha^a)$ by
$$M^a:=(M\oplus\SO^*_X)/P,\quad \alpha^a(x,u)=u\cdot\alpha(x)$$
where $P=\{(x,\alpha(x)^{-1})\,|x\in\alpha^{-1}(\SO^*_X)\}$. Let $f:X\to Y$
be a morphism of schemes. For a log structure $M$ on $X$, we can define a 
log structure on $Y$ called the direct image of $M$, to be the fibre product 
of sheaves   
$$\CD   
 @.  f_*M \\
  @.    @VVV     \\
 \SO_Y @>>> f_*\SO_X
\endCD.$$
For a log structure $N$ on $Y$, we define a log structure $f^*N$ on $X$ called
the inverse image of $N$ to be the log structure associated to the pre-log
structure $$f^{-1}(N)\to f^{-1}(\SO_Y)\to \SO_X.$$
 
\proclaim{Definition 1.1} Let $\alpha: M\to\SO_X$ and $\beta:N\to\SO_Y$ be
pre-log structures and $f:(X,M)\to (Y,N)$ be a morphism of log schemes. Then
the $\SO_X$-module $\Omega^1_{X/Y}(log(M/N))$ called logarithmic differential
sheaf is defined to be the quotient of
$$\Omega^1_{X/Y}\oplus(\SO_X\otimes_{\Bbb Z}M^{gp})$$
($\Omega^1_{X/Y}$ is the usual relative differential module) divided by the 
$\SO_X$-submodule generated locally by local sections of the following forms
\roster\item $(d\alpha(a),0)-(0,\alpha(a)\otimes a)$ with $a\in M$.
\item $(0,1\otimes a)$ with $a\in\text{Image$(f^{-1}(N)@>h>>M)$}$.\endroster
\endproclaim

It is easily seen that if $M^a$ and $N^a$ denote the associated log structures
respectively, we have
$$\Omega^1_{X/Y}(log(M/N))=\Omega^1_{X/Y}(log(M^a/N))=\Omega^1_{X/Y}(log(M^a/N^a)).$$
We collect some easy facts in the following proposition which may be useful in
the paper.

\proclaim{Proposition 1.1} Let $f:X\to Y$ be a morphism of schemes, and $N^a$
the log structure associated to a pre-log structure $N$ on $Y$. Then 
\roster\item $f^*(N^a)$ coincides with the log structure associated to the
pre-log structure $f^{-1}(N)\to \SO_X.$
\item If $f:(X,M^a)\to (Y,N^a)$ is a morphism of log schemes such that
$$f^{-1}(N)\to M$$
is surjective, we have 
$\Omega^1_{X/Y}(log(M/N))=\Omega^1_{X/Y}$ (the usual relative differential sheaf).
\item If we have a cartesian diagram of log schemes
$$\CD
(X',M') @>f>>(X,M)\\
@VVV          @VVV \\
(Y',N')@>>> (Y,N),  
\endCD$$
we have an isomorphism $f^*\Omega^1_{X/Y}(log(M/N))\cong\Omega^1_{X'/Y'}(log(M'/N'))$.
\endroster\endproclaim

Now we are going to introduce the sheaf of differential operators on general log
schemes although we need it only for some special log structures, we hope this
general treatment to be useful in the future. Fix a morphism $X^{\dagger}
=(X,M)\to Y^{\dagger}=(Y,N)$
of log schemes and denote the sheaf of $\SO_Y$-derivations of $\SO_X$ by $T_{X/Y}$. 
We write $\Omega^1_{X/Y}(log)$ simply for
$\Omega^1_{X/Y}(log(M/N))$, and 
$$\bar d:\SO_X @>d>>\Omega^1_{X/Y}\to\Omega^1_{X/Y}(log)$$
denotes the canonical logarithmic derivation. 

\proclaim{Definition 1.2} A derivation $\delta\in T_{X/Y}$ is called a logarithmic 
derivation if there exists a 
$\theta\in Hom_{\SO_X}(\Omega^1_{X/Y}(log),\SO_X)$ such that   
$$\CD
\SO_X @>\delta>>\SO_X \\
@V{\bar d}VV     @A{\theta}AA \\
\Omega^1_{X/Y}(log)@=\Omega^1_{X/Y}(log)
\endCD$$
is commutative.\endproclaim

\remark{Remark 1.1} The sheaf of logarithmic derivations, 
denoted by $T_{X/Y}(log)$, is a subsheaf of $T_{X/Y}$. By definition, we have
a surjection
$$\aligned\Cal Hom_{\SO_X}(\Omega^1_{X/Y}(log),\SO_X)&\to T_{X/Y}(log)\\
u&\mapsto u\circ\bar d,\endaligned$$
which is not injective in general since $\Omega^1_{X/Y}(log)$ is not
generated by $\{\bar df\}_{f\in\SO_X}$. However, we will see in Lemma 1.2
that for all the logarithmic structures we concern in this paper the above 
surjection is actually an isomorphism.\endremark 

Let $\Cal End_{\SO_Y}(\SO_X)$ be sheaf of
$\SO_Y$-linear maps and $\SO_X\subset \Cal End_{\SO_Y}(\SO_X)$ be the
subsheaf of maps multiplying by elements of $\SO_X$. It is clear that
$\Cal End_{\SO_Y}(\SO_X)$ is a sheaf of noncommutative rings, thus it
has two $\SO_X$-module structures (left and right multipications).

\proclaim{Definition 1.3} The sheaf 
$\SD_{X^{\dagger}/Y^{\dagger}}\subset\Cal End_{\SO_Y}(\SO_X)$ 
of subrings generated by $\SO_X$ and $T_{X/Y}(log)$
is called the sheaf of logarithmic differential operators on log scheme
$X^{\dagger}=(X,M)$ over $Y^{\dagger}=(Y, N)$. We will simply call it the sheaf of 
log differential
operators. For any integer $k>0$, we define inductively the sheaf
of k-th log differential operators:
$$\SD^0_{X^{\dagger}/Y^{\dagger}}:=\SO_X,$$
$$\SD^k_{X^{\dagger}/Y^{\dagger}}:=
\{D\in\SD_{X^{\dagger}/Y^{\dagger}}\,|\,[D,f]\in\SD^{k-1}_{X^{\dagger}/Y^{\dagger}}
,\,\text{for any $f\in\SO_X$}\}.$$ \endproclaim  

If $\Omega^1_{X/Y}(log)$ is locally free, we can describe $\SD_{X^{\dagger}/Y^{\dagger}}$
locally. Let $$dlog:\SO_X\otimes_{\Bbb Z}M^{gp}\to\Omega^1_{X/Y}(log)$$
be the surjection and choose locally $t_1,...,t_r\in M$ such that 
$\{dlog(t_i)\}_{1\le i\le r}$ is an $\SO_X$-base of $\Omega^1_{X/Y}(log)$.
Let $f_i=\alpha(t_i)\in\SO_X$ ($i=1,...,r$), where $\alpha:M\to\SO_X$ is the
log structure. Then there exists locally a system of generators $\{\p_i\}_{1\le i\le r}$
of $T_{X/Y}(log)$ such that
$$\p_i(f_j)=\cases f_i&\quad j=i,\\ 0&\quad j\neq i.\endcases$$
All $[\p_i,\p_j]$ vanish on the subring $\SO'\subset\SO_X$ generated by
$f_1,...,f_r$ and $\SO_Y$, which means that $[\p_i,\p_j]$ are derivations of
$\SO_X$ over $\SO'$, thus $[\p_i,\p_j]$ vanish on $\SO_X$ since 
$\Omega^1_{\SO_X/\SO'}$ is a torsion sheaf. Therefore any local section
$D\in\SD_{X^{\dagger}/Y^{\dagger}}$ can be expressed as a finite sum
$$D=\sum\lambda_{\beta_1,...,\beta_r}\p_1^{\beta_1}\cdots\p_r^{\beta_r}.\tag1.1$$
We introduce a notation $[D,a_1\star\cdots\star a_n]$
for any local section $D\in\SD_{X^{\dagger}/Y^{\dagger}}$ and $a_1,...,a_n\in \SO_X$. 
The $[D,a_1\star\cdots\star a_n]\in \SD_{X^{\dagger}/Y^{\dagger}}$ is defined inductively by
$$[D,a_1\star\cdots\star a_n]:=[[D,a_1\star\cdots\star a_{n-1}],a_n].$$
If $a_1=\cdots=a_n$, we write $[D,a_1\star\cdots\star a_n]=[D,a_1^{\star n}]$,
thus the notation 
$[D,a_1^{\star i_1}\star a_2^{\star i_2}\star\cdots\star a_k^{\star i_k}]$ is clear. 
Notice that $[\p_i,f_j]=0$ for $i\neq j$, one checks easily that
$$[\p_1^{\beta_1}\cdots\p_r^{\beta_r},f_i^{\star\beta_i}]=
\p_1^{\beta_1}\cdots\p_{i-1}^{\beta_{i-1}}\cdot[\p_i^{\beta_i},
f_i^{\star\beta_i}]\cdot\p_{i+1}^{\beta_{i+1}}\cdots\p_r^{\beta_r}.$$
Thus
$[\p_1^{\beta_1}\cdots\p_r^{\beta_r},f_1^{\star\beta_1}\star\cdots\star f_r^{\star\beta_r}]=
[\p_1^{\beta_1},f_1^{\star\beta_1}]\cdots[\p_r^{\beta_r},
f_r^{\star\beta_r}].$ One computes that $[\p_i^{\beta_i},f_i^{\star\beta_i}]=
(\beta_i!)f_i^{\beta_i}$ and $[\p_i^{k_i},f_i^{\star\beta_i}]=0$ if $k_i<\beta_i.$  
For any $D\in\SD_{X^{\dagger}/Y^{\dagger}}^k,$ 
let $\lambda_{\alpha_1,...,\alpha_r}\p_1^{\alpha_1}\cdots\p_r^{\alpha_r}$ be the
first term of $D$ in (1.1) (namely $\alpha_1\ge\beta_1,...,\alpha_r\ge\beta_r$ and 
at least one inequality is strict). Then
$$\aligned [D,f_1^{\star\alpha_1}\star\cdots\star f_r^{\star\alpha_r}]&=
\lambda_{\alpha_1,...,\alpha_r}[\p_1^{\alpha_1}\cdots\p_r^{\alpha_r},
f_1^{\star\alpha_1}\star\cdots\star f_r^{\star\alpha_r}]\\&=
(\alpha_1!)\cdots(\alpha_r!)\lambda_{\alpha_1,...,\alpha_r}f_1^{\alpha_1}\cdots
f_r^{\alpha_r}.\endaligned$$
If $\alpha_1+\cdots+\alpha_r>k$, then 
$[D,f_1^{\star\alpha_1}\star\cdots\star f_r^{\star\alpha_r}]=0$ by definition of 
$\SD^k_{X^{\dagger}/Y^{\dagger}}$, thus
$$\lambda_{\alpha_1,...,\alpha_r}\cdot f_1^{\alpha_1}\cdots
f_r^{\alpha_r}=0.\tag1.2$$

\proclaim{Proposition 1.2} Let 
$Gr_i(\SD_{X^{\dagger}/Y^{\dagger}})=
\SD^i_{X^{\dagger}/Y^{\dagger}}/\SD^{i-1}_{X^{\dagger}/Y^{\dagger}}$ and 
$T(Gr_1(\SD_{X^{\dagger}/Y^{\dagger}}))$
be the tensor algebra. Then\roster
\item The symbol map $\sigma: Gr_1(\SD_{X^{\dagger}/Y^{\dagger}})\to T_{X/Y}(log),$
defined by $\sigma(D)(f)=[D,f](1)$, is an isomorphism.
\item $\SD_{X^{\dagger}/Y^{\dagger}}=\bigcup^{\infty}_{i=0}
\SD^i_{X^{\dagger}/Y^{\dagger}}$, $\SD^i_{X^{\dagger}/Y^{\dagger}}
\cdot\SD^j_{X^{\dagger}/Y^{\dagger}}
\subset\SD^{i+j}_{X^{\dagger}/Y^{\dagger}}$ and $f^{-1}(\SO_Y)$ is in
the center of $\SD_{X^{\dagger}/Y^{\dagger}}$.
\item The left and right $\SO_X$-module structures on 
$Gr_i(\SD_{X^{\dagger}/Y^{\dagger}})$
coincide 
\item The natural map 
$$T(Gr_1(\SD_{X^{\dagger}/Y^{\dagger}}) )
\to Gr(\SD_{X^{\dagger}/Y^{\dagger}}):=
\bigoplus^{\infty}_{i=0}Gr_i(\SD_{X^{\dagger}/Y^{\dagger}})$$ is surjective.
\endroster\endproclaim

\demo{Proof} One checks that $\sigma(D)$ is a derivation of $\SO_X$. To see it
factorizing through $\bar d:\SO_X\to\Omega^1_{X/Y}(log)$, it is enough to show
that for any $m\in M$ there exists a unique $\theta(dlog(m))\in \SO_X$ such that
$\sigma(D)(\alpha(m))=\alpha(m)\theta(dlog(m)),$ which can be checked by using
the fact that $D$ is a map composed by logarithmic derivations.
 
The proof of others is easy, we will omit it but just remark that we
are not able to claim the natural map in (4) induces surjections
$$T^i(Gr_1(\SD_{X^{\dagger}/Y^{\dagger}}) )\to Gr_i(\SD_{X^{\dagger}/Y^{\dagger}}).$$ 
\enddemo

A reduced scheme $Y$ is called a normal crossing variety of dimension $d$
if the completion of local ring $\SO_{Y,y}$ at each point $y$ is isomorphic to
$\Bbb C\{x_0,...,x_d\}/(x_0\cdots x_{r})$ for some $r=r(y)$ such that 
$0\le r\le d.$ Now we discuss the log structures on normal crossing varieties
and on smooth varieties induced by normal crossing divisors. These are the 
only log structures we will concern in this paper. 

\proclaim{Lemma 1.1} Let $f:X\to S$ be a flat family of (reduced) normal 
crossing varieties of dimension $d$, and $X$, $S$ smooth. Assume that 
$\Delta\subset S$ is a normal crossing divisor and $Y:=f^{-1}(\Delta)\subset X$ 
such that $f:X\ssm Y\to S\ssm\Delta$ is smooth. Then, for any $x\in X$, we can choose
isomorphisms $$\hat\SO_{X,x}\cong\Bbb C\{x_1,...,x_{d+1},...,x_{d+m}\},\quad
\hat\SO_{S,f(x)}\cong\Bbb C\{\pi_1,...,\pi_m\}$$ 
such that
$$f^{\sharp}:\Bbb C\{\pi_1,...,\pi_m\}\to\Bbb C\{x_1,...,x_{d+1},...,x_{d+m}\}$$
is a $\Bbb C$-algebra homomorphism with 
$$f^{\sharp}(\pi_2)=x_{d+2},\quad.\,\,.\,\,.\quad, f^{\sharp}(\pi_m)=x_{d+m}$$
and the local equation of $\Delta$ at $f(x)$ is
$$\pi_{i_1}\cdot\pi_{i_2}\cdots\pi_{i_s}=0.$$\endproclaim

\demo{Proof} Let $\hat\SO_{S,f(x)}\cong\Bbb C\{\pi_1,...,\pi_m\}$ and
$\hat\SO_{X,x}\cong\Bbb C\{z_1,...,z_{d+m}\}$. Then, by definition,
$$\varphi:\frac{\Bbb C\{y_1,...,y_{d+1}\}}{(y_1\cdots y_{d+1})}\cong
\frac{\Bbb C\{z_1,...,z_{d+m}\}}{(f^{\sharp}(\pi_1),...,f^{\sharp}(\pi_m))}.$$
We can assume that $r>1$ (otherwise, $f$ is smooth at $x$, the lemma is clear).
Let
$$\varphi(\bar y_i)=\overline{\varphi_i(z_1,...,z_{d+m})}\in
\frac{\Bbb C\{z_1,...,z_{d+m}\}}{(f^{\sharp}(\pi_1),...,f^{\sharp}(\pi_m))},$$
we can write $\varphi_i(z_1,...,z_{d+m})=\sum_{j=1}^{d+m}a_{ij}z_j\,+\,
\varphi_i^{\ge 2}$, where $\varphi_i^{\ge 2}$ denotes the part of 
$\varphi_i(z_1,...,z_{d+m})$ with order $\ge 2$.
Then $$\frac{\p(\varphi_1,...,\varphi_{d+1})}{\p(z_1,...,z_{d+m})}\mid_{(0,...,0)}
=\pmatrix
a_{11}&\hdots&a_{1\,d+m}\\
\vdots&\ddots&\vdots\\
a_{d+1\,1}&\hdots&a_{d+1\,d+m}
\endpmatrix$$
has rank $d+1$. Otherwise, there is $0\neq (k_1,...,k_{d+1})\in\Bbb C^{d+1}$ 
such that
$$\sum_{i=1}^{d+1}k_i\varphi_i(z_1,...,z_{d+m})=\sum_{i=1}^{d+1}k_i\varphi_i^
{\ge 2}\in(z_1,...,z_{d+m})^2,$$
which implies that
$$\sum^{d+1}_{i=1}k_i\bar y_i=
\varphi^{-1}(\overline{\sum_{i=1}^{d+1}k_i\varphi_i(z_1,...,z_{d+m})})\in
(\bar y_1,...,\bar y_{d+1})^2.$$
Thus there is $g(y_1,...,y_{d+1})\in (y_1,...,y_{d+1})^2$ such that
$$\sum^{d+1}_{i=1}k_i y_i\,-\,g(y_1,...,y_{d+1})\in (y_1\cdots y_r),$$
which is impossible since $r>1$. Replaceing some $z_{j_1},...,z_{j_{d+1}}$ by
$\varphi_1(z_1,...,z_{d+m})$, ..., $\varphi_{d+1}(z_1,...,z_{d+m})$, we can assume
that $$\varphi:\frac{\Bbb C\{y_1,...,y_{d+1}\}}{(y_1\cdots y_{d+1})}\cong
\frac{\Bbb C\{z_1,...,z_{d+m}\}}{(f^{\sharp}(\pi_1),...,f^{\sharp}(\pi_m))}$$
such that $\varphi(\bar y_i)=\bar z_i$ ($i=1,...,d+1$). Let
$$\varphi^{-1}(\bar z_{d+j})=\overline{g_j(y_1,...,y_{d+1})}\in
\frac{\Bbb C\{y_1,...,y_{d+1}\}}{(y_1\cdots y_r)},$$
then $z_{d+j}-g_j(z_1,...,z_{d+1})\in (f^{\sharp}(\pi_1),...,f^{\sharp}(\pi_m)).$
If we write that (for $j=2,...,m$)
$$z_{d+j}-g_j(z_1,...,z_{d+1})=\sum^m_{i=1}a_{ji}f^{\sharp}(\pi_i)\,+\,
\text{higher order terms},$$
then 
$$\pmatrix
1&0&\hdots&0\\
0&1&\hdots&0\\
\vdots&\vdots&\ddots&\vdots\\
0&0&\hdots&1\endpmatrix=
\pmatrix
a_{21}&\hdots&a_{2\,m}\\
\vdots&\ddots&\vdots\\
a_{m\,1}&\hdots&a_{m\,m}
\endpmatrix \cdot
\frac{\p(f^{\sharp}(\pi_1),...,f^{\sharp}(\pi_m))}{\p(z_{d+2},...,z_{d+m})}
\mid_{(0,...,0)}.$$
 Thus
$$\frac{\p(f^{\sharp}(\pi_1),...,f^{\sharp}(\pi_m))}{\p(z_{d+2},...,z_{d+m})}
\mid_{(0,...,0)}$$
has rank $m-1$, and we can choose the isomorphism
$$\hat\SO_{X,x}\cong\Bbb C\{x_1,...,x_{d+1},...,x_{d+m}\}$$
such that $f^{\sharp}(\pi_2)=x_{d+2}$, ..., $f^{\sharp}(\pi_m)=x_{d+m}$ by 
changing the order of $\pi_1$, ..., $\pi_m$.\enddemo

We are not able to prove in Lemma 1.1 that $Y=f^{-1}(\Delta)$ is
a normal crossing divisor in $X$ even if each fibre is a normal crossing
variety. Assume that $Y$ is a normal crossing divisor, we have canonical
log structures on $X$ and $S$ 
$$logY=\{g\in\SO_X\,|\,\text{$g$ is invertible outside $Y$}\}\subset\SO_X$$     
$$log\Delta=\{g\in\SO_S\,|\,\text{$g$ is invertible outside $\Delta$}\}\subset\SO_S.$$
These are fine log structures, and if one writes locally that
$$Y=\bigcup^r_{i=1}\{x_i=0\}\quad\Delta=\bigcup^s_{i=1}\{\pi_i=0\},$$
then $logY$ and $log\Delta$ are associated to the pre-log. structures ([KK]):
$$\aligned \Bbb N^r&\to\SO_X\\(n_i)_{1\le i\le r}&\mapsto \prod x_i^{n_i},
\endaligned \qquad \aligned \Bbb N^s&\to\SO_S\\(n_i)_{1\le i\le s}
&\mapsto \prod \pi_i^{n_i}.\endaligned$$
For the local descriptions and properties of $\Omega^1_S(log\Delta)$ and
$\Omega^1_X(logY)$, we refer to [EV1] and [EV2]. We will use 
$\Omega^1_{X/S}(logY)$ to denote $\Omega^1_{X/S}(logY/\Delta)$. 

\proclaim{Proposition 1.3} Let $f:X\to S$ be a flat family of normal
crossing varieties of dimension $d$, $X$ and $S$ be smooth. Assume that
$Y:=f^{-1}(\Delta)$ is a normal crossing divisor such that
$$f:X\ssm Y\to S\ssm\Delta$$
is smooth. Then we have the associated exact sequence
$$0\to f^*\Omega^1_S(log\Delta)@>j>>\Omega^1_X(logY)\to
\Omega^1_{X/S}(logY)\to 0,$$
and the following are equivalent\roster
\item $(X,logY)@>f>> (S,log\Delta)$ is log smooth.
\item The image of $j$ is locally a direct summand.
\item  For any singular point $x\in X$ of $f$, we can choose coordinates
$$ 
\hat\SO_{S,f(x)}\cong\Bbb C\{\pi_1,...,\pi_m\}\hookrightarrow 
\Bbb C\{x_1,...,x_{d+1},\pi_2,...,\pi_m\}\cong\hat\SO_{X,x}
$$
such that $\pi_1=x_1\cdots x_r$ for some $1\le r\le d+1$ and
$$\pi_1\cdot\pi_2\cdots\pi_s=0$$ is the local equation of $\Delta$ at $f(x)$.
\item $\Omega^1_{X/S}(logY)$ is locally free.\endroster
\endproclaim

\demo{Proof} The map $j$ has to be injective since it is injective at the 
generic point of $S$ and $\Omega^1_{X/S}(log\Delta)$ is locally free.

The $(1)\Leftrightarrow (2)$ follows the Proposition (3.12) of [KK], and
$(2)\Leftrightarrow (4)$ is obvious. We prove that $(3)\Leftrightarrow(4).$
Firstly, $(3)\Rightarrow (4)$ is clear since $\Omega^1_{X/S}(logY)$
is locally isomorphic to
$$\frac{\hat\SO_{X,x}\{dx_1,...,dx_{d+1},e_1,...,e_r\}}
{(dx_1-x_1e_1,...,dx_r-x_re_r,e_1+\cdots+e_r)},$$
which is a free module generated by 
$\{\frac{1}{x_1}dx_1,...,\frac{1}{x_{r-1}}dx_{r-1},dx_{r+1},...,dx_{d+1}\}.$
To prove $(4)\Rightarrow (3)$, we only need to show that local equation of
$\Delta$ is divisible by $\pi_1$ (assume that we choose the coordinates as in 
Lemma 1.1). If it is not so, we can assume the local equation of $\Delta$ to be
$$\pi_2\cdots\pi_s=0.$$
Then, as in the Lemma 1.1, $f^{\sharp}(\pi_2\cdots\pi_s)=x_{d+2}\cdots x_{d+s}$
is the local equation of $Y$ at $x$, and $\Omega^1_{X/S}(logY)$ is locally isomorphisc to
$$\frac{\hat\SO_{X,x}\{dx_1,...,dx_{d+1}\}}{\sum^{d+1}_{i=1}
\frac{\p f^{\sharp}(\pi_1)}{\p x_i}dx_i}.$$
It is not locally free except one of $\frac{\p f^{\sharp}(\pi_1)}{\p x_i}$
is invertible, which means that $f$ is smooth at $x$. But (3) is clear 
at this case.\enddemo

Let $X_s$ ($s\in\Delta$) be a fibre of $f:X\to S$, then $X_s$ is a normal
crossing variety with log structure $M:=(logY)|_{X_s}@>\alpha>>\SO_{X_s}$ such
that $X_s^{\dagger}:=(X_s,M)$ is a log smooth variety. We will describe $M$ locally and 
to show how it gives a log structure in the sense of [KN].

For $x\in X_s$ a singular point of $X_s$, there is a neighbourhood 
$U_{\lambda}\subset X$ of $x$ and holomorphic functions 
$x^{\lambda}_i\in\SO_X(U_{\lambda})$ such that
$$\aligned &U_{\lambda}\to\Bbb C^{d+m}\\
&p\mapsto (x_1^{\lambda}(p),...,x_{d+m}^{\lambda})\endaligned$$
is an open embedding and $X_s\cap U_{\lambda}\subset U_{\lambda}$ is defined
by $$x_1^{\lambda}\cdots x_r^{\lambda}=0,\quad x_{d+2}^{\lambda}=\cdots=
x_{d+m}^{\lambda}=0,$$
$Y\cap U_{\lambda}$ is defined by 
$x_1^{\lambda}\cdots x_r^{\lambda}\cdot x_{d+i_1}^{\lambda}\cdots x_{d+i_s}^{\lambda}=0,$
for some $1<r\le d+1$ and $2\le i_1<\cdots<i_s\le m$. The log structure
$logY|_{U_{\lambda}}$ on $U_{\lambda}$ is associated to the pre-log structure
$$\aligned  \Bbb N^{r+s}&\to\SO_{U_{\lambda}}\\
(n_1,...,n_r,n_{d+i_1},...,n_{d+i_s})&\mapsto 
(x_1^{\lambda})^{n_1}\cdots (x_r^{\lambda})^{n_r}\cdot 
(x_{d+i_1}^{\lambda})^{n_{d+i_1}}\cdots (x_{d+i_s}^{\lambda})^{n_{d+i_s}}.\endaligned$$

Let $z_i^{\lambda}=x_i^{\lambda}|_{X_s\cap U_{\lambda}}
\in\SO_{X_s}(X_s\cap U_{\lambda})$, then $M|_{X_s\cap U_{\lambda}}$ is associated 
the pre-log structure
$$\aligned \Bbb N^r\oplus\Bbb N^s=\Bbb N^{r+s}
&@>\upsilon>>\SO_{U_{\lambda}\cap X_s}\\
(n_1,...,n_r,n_{d+i_1},...,n_{d+i_s})&\mapsto (z_1^{\lambda})^{n_1}
\cdots (z_r^{\lambda})^{n_r}\cdot (z_{d+i_1}^{\lambda})^{n_{d+i_1}}
\cdots (z_{d+i_s}^{\lambda})^{n_{d+i_s}}.\endaligned$$
Clearly $ker(\upsilon)=\Bbb N^s$ and $\upsilon^*(\SO_{U_{\lambda}\cap X_s}^*)=
\{(0,...,0)\in\Bbb N^{r+s}\}$, so
$M|_{X_s\cap U_{\lambda}}=\Bbb N^{r+s}\oplus\SO^*_{U_{\lambda}\cap X_s}$ and
$\alpha_{\lambda}:M|_{X_s\cap U_{\lambda}}\to\SO_{U_{\lambda}\cap X_s}$ is 
defined by $\alpha_{\lambda}(\vec n, g)=\upsilon(\vec n)\cdot g.$
Thus we get a partial open cpvering $\{V_{\lambda}:=U_{\lambda}\cap X_s\}$ of
$X_s$ containing the singular locus of $X_s$ and systems of holomorphic 
functions 
$$z_i^{(\lambda)}:=\alpha_{\lambda}(e_i,1)\in \SO_{X_s}(V_{\lambda}),$$
where $i=1,...,r_{\lambda}$ and $e_i=(0,...,1,...,0)\in\Bbb N^{r_{\lambda}+s}.$
These functions satisfy\roster
\item There is an isomorphism $\varphi_{\lambda}$ from $V_{\lambda}$ to an open
neighborhood of $(0,...,0)$ of the variety
$$\{(x_1,...,x_{d+1})\in\Bbb C^{d+1}\,|\,x_1\cdots x_{r_{\lambda}}=0\}$$
such that $\varphi^*_{\lambda}(x_j)=z_j^{(\lambda)}$ for $1\le j\le r_{\lambda}$.
\item If $V_{\lambda}\cap V_{\mu}\neq\empty$ ($r_{\lambda}=r_{\mu}$ at this case), 
then there exist
invertible holomorphic functions $u_j^{(\lambda\mu)}$ ($1\le j\le r_{\lambda}$)
on $V_{\lambda}\cap V_{\mu}$ and a permutation $\sigma\in S_{r_{\lambda}}$
such that 
$$\text{$z^{(\lambda)}_{\sigma(j)}=u_j^{(\lambda\mu)}z_j^{(\mu)}$ and 
$u_1^{(\lambda\mu)}\cdots u_{r_{\lambda}}^{(\lambda\mu)}=1$ on $V_{\lambda}\cap
V_{\mu}$}.$$\endroster
Thus we have a $log$ $atlas$ in the sense of [KN]. To check (2), noting that
$$\aligned Y&=\bigcup^{r_{\lambda}}_{i=1}\{x_i^{\lambda}=0\}\cup
\bigcup^s_{k=1}\{x^{\lambda}_{d+i_k}=0\}\\
&=\bigcup^{r_{\lambda}}_{i=1}\{x_i^{\mu}=0\}\cup
\bigcup^s_{k=1}\{x^{\mu}_{d+i_k}=0\}\endaligned$$
on $U_{\lambda}\cap U_{\mu}$, we have $\sigma\in S_{r_{\lambda}+s}$ and invertible 
holomorphic functions $u_j^{(\lambda\mu)}\in\SO_X(U_{\lambda}\cap 
U_{\mu})$ such that $$x^{\lambda}_{\sigma(j)}=u_j^{(\lambda\mu)}x_j^{\mu}.$$
Since $x^{\lambda}_{d+i_k}=x^{\mu}_{d+i_k}=0$ on 
$X_s\cap U_{\lambda}\cap U_{\mu}$, we have $\sigma(j)\in\{1,...,r_{\lambda}\}$ 
if $j\in\{1,...,r_{\lambda}\}.$
On the other hand, for fixed $(\pi_1,...,\pi_m)=m_s$, we have
$\pi_1=x^{\lambda}_1\cdots x^{\lambda}_{r_{\lambda}}=x^{\mu}_1\cdots 
x^{\mu}_{r_{\lambda}},$
thus $u_1^{(\lambda\mu)}\cdots u_{r_{\lambda}}^{(\lambda\mu)}=1$. The restrictions of 
$\{u_j^{(\lambda\mu)}\}$ to $X_s\cap U_{\lambda}\cap U_{\mu}$
give the required functions in (2).

\proclaim{Lemma 1.2} For the smooth variety $X$ with log structure $logY$ given
by a normal crossing divisor $Y$, and for the normal crossing variety $X_s$ with
log structure in the sense of [KN] (in particular, for $logY|_{X_s}$), the
surjections
$$\aligned\Omega^1_{X/S}(logY)^*&\to T_{X/S}(logY)\\
u&\mapsto u\circ\bar d\endaligned
\quad\aligned\Omega^1_{X_s}(log)^*&\to T_{X_s}(log)\\
u&\mapsto u\circ\bar d\endaligned$$
are isomorphisms, where $\SE^*$ denotes the dual of $\SE$.
\endproclaim

\demo{Proof} It is enough to check the injectivity of the above morphisms, which
is a local problem. Thus we may assume that
$$\SO_{X_s}=\frac{\Bbb C[[x_1,...,x_r,...,x_n]]}{(x_1\cdots x_r)},\quad
\Omega^1_ {X_s}(log)=\frac{\SO_{X_s}\{\frac{\bar dx_1}{x_1},...,
\frac{\bar dx_r}{x_r},\bar dx_{r+1},...,\bar dx_n\}}{(\frac{\bar dx_1}{x_1}+\cdots +
\frac{\bar dx_r}{x_r})}.$$
For any $u:\Omega^1_{X_s}(log)\to\SO_{X_s}$ such that $u(\bar df)=0$ for
all $f\in\SO_{X_s}$, we need to show that for $i=1,...,r$
$$\bar a_i:=u(\frac{1}{x_i}\bar dx_i)=0.$$
Firstly, $a_i=b_i(x_1\cdots\hat x_i\cdots x_r)$ for some $b_i\in\Bbb C[[x_1,...,x_n]]$ 
since $\bar a_i\bar x_i=u(\bar dx_i)=0$. Secondly,
$a_1+\cdots+a_r$ is divisible by $x_1\cdots x_r$ since
$$\bar a_1+\cdots+\bar a_r=u(\frac{1}{x_1}\bar dx_1+\cdots+\frac{1}{x_r}\bar dx_r)=0.$$
These facts imply that $b_i$ is dvisible by $x_i$ ($i=1,...,r$), which means
 all $\bar a_i$ are zero.

When $X$ is a smooth variety, local ring (analytic) at any point of $X$ is
integral, the injectivity is obvious from the above proof. 
\enddemo 

Thus for a logarithmic derivation $\theta$, we will write $<\theta,\cdot>$
denoting the unique element of $\Omega^1_{X_s}(log)^*$ such that 
$<\theta,\bar da>=\theta(a)$.

\proclaim{Proposition 1.4} Let $W\subset X_s$ be the singular locus of $X_s$
and $\Cal I_W\subset\SO_{X_s}$ the ideal sheaf of $W$. Then
every local section $D\in \SD^k_{X_s^{\dagger}}$ can be expressed
uniquely into $$D=\sum_{\beta_1+\cdots+\beta_d\le k}
\lambda_{\beta_1,...,\beta_d}\p_1^{\beta_1}\cdots\p_d^{\beta_d}$$
for a local basis $\p_1$, ..., $\p_d$ of $T_{X_s}(log)$ and 
$D(\Cal I_W)\subset I_W.$
In particular, we have the canonical exact sequence
$$0\to \SD^{k-1}_{X_s^{\dagger}}\to \SD^k_{X_s^{\dagger}}
@>\sigma_k>> S^kT_{X_s}(log)\to 0,\tag1.3$$
where $S^kT_{X_s}(log)$ is the subsheaf of symmetric tensors of $T^k(T_{X_s}(log))$ 
and $\sigma_k$ will be defined in the proof.              
\endproclaim

\demo{Proof}  Locally, $\hat\SO_{X_s,x}=\Bbb C\{x_1,...,x_{d+1}\}/(x_1\cdots x_r)$, 
and $T_{X_s}(log)$ is locally generated by
$$x_1\frac{\p\,\,}{\p x_1},...,x_r\frac{\p\,\,}{\p x_r},
\frac{\p\,\,\,}{\p x_{r+1}},...,\frac{\p\,\,\,}{\p x_{d+1}}$$
with a relation $$x_1\frac{\p\,\,}{\p x_1}+\cdots +x_r\frac{\p\,\,}{\p x_r}=0.$$
If we take the local basis $\p_1$, ..., $\p_d$ of $T_{X_s}(log)$ to be
$$\p_1=x_1\frac{\p\,\,}{\p x_1},...,\p_{r-1}=x_{r-1}\frac{\p\,\,\,}{\p x_{r-1}},
\p_r=\frac{\p\,\,\,}{\p x_{r+1}},...,\p_d=\frac{\p\,\,\,}{\p x_{d+1}},$$
one can check that $\p_i^{\alpha_i}=x_ig_i(x_i,\frac{\p\,\,}{\p x_i})$
when $\alpha_i>0$ and $1\le i<r$, where $g_i(x,y)$ is a polynoimal. Thus the first term of
(1.1) becomes into 
$$\lambda_{\alpha_1,...,\alpha_d}\p_1^{\alpha_1}\cdots\p^{\alpha_d}_d=  
\lambda_{\alpha_1,...,\alpha_d}x_{i_1}\cdots x_{i_t}\cdot g_{i_1}\cdots 
g_{i_t}\cdot\p_r^{\alpha_r}\cdots\p_d^{\alpha_d},$$
and the equality (1.2) becomes into 
$\lambda_{\alpha_1,...,\alpha_d}x^{\alpha_{i_1}}_{i_1}\cdots x^{\alpha_{i_t}}_{i_t}=0,$ 
where $\alpha_{i_1},...,\alpha_{i_t}$ are nonzero
integers of $\alpha_1,...,\alpha_{r-1}$. This implies that 
$\lambda_{\alpha_1,...,\alpha_d}$ is divisible by all $x_i$ with 
$i\in\{1,...,r\}\ssm\{i_1,...,i_t\}$, namely 
$\lambda_{\alpha_1,...,\alpha_d}\p_1^{\alpha_1}\cdots\p^{\alpha_d}_d=0,$ which
means that any $D\in \SD^k_{X_s^{\dagger}}$ can be expressed
uniquely into $$D=\sum_{\beta_1+\cdots+\beta_d\le k}
\lambda_{\beta_1,...,\beta_d}\p_1^{\beta_1}\cdots\p_d^{\beta_d}.$$

On the other hand, $\Cal I_W$ is locally generated by $\frac{x_1\cdots x_r}{x_i}$ 
($i=1,...,r$), and
$$\p_j(\frac{x_1\cdots x_r}{x_i})=\cases 0, & \text{if $j=i$ or $j\ge r$}\\
\frac{x_1\cdots x_r}{x_i}. & \text{if $i\neq j<r$}\endcases$$
Thus $D(\Cal I_W)\subset \Cal I_W$. Now we see that the natural map in
Proposition 1.2 (4) induces isomorphisms
$$S^k(Gr_1(\SD_{X^{\dagger}/Y^{\dagger}}))\to Gr_k(\SD_{X^{\dagger}/Y^{\dagger}}).$$
For any local section $D\in \SD^k_{X_s^{\dagger}}$, let $D_{<k}$ denote
the lower order part and write
$$D=D_{<k}+\sum_{\beta_1+\cdots+\beta_d= k}
\lambda_{\beta_1,...,\beta_d}\p_1^{\beta_1}\cdots\p_d^{\beta_d},$$
$$D(\omega_1,...,\omega_k)=\sum_{\beta_1+\cdots+\beta_d= k}
\lambda_{\beta_1,...,\beta_d}
\left(\prod^{\beta_1}_1<\sigma(\p_1),\omega_i>
\cdots\prod_{\beta_{d-1}+1}^{\beta_d}<\sigma(\p_d),\omega_i>\right)$$
where $\omega_1,...,\omega_k $ are elements of $\Omega^1_{X_s}(log).$  The symbol 
$\sigma_k(D)$ of $D$ as a symmetric function on $\otimes^k\Omega^1_{X_s}(log)$
is defined to be
$$\sigma_k(D)(\omega_1,...,\omega_k)=
\sum_{\tau\in S_k}D(\omega_{\tau(1)},...,\omega_{\tau(k)}).$$
This gives the exact sequence (1.3) and coincides the definition of [GJ] and
[We] in smooth case (see Remark 2.2.4 of [GJ]). 
\enddemo

\heading \S2 Logarithmic heat operators and logarithmic connections\endheading

In this section, we generalize the definitions and arguments about heat
operators and connections in [GJ] and [We] to the logarithmic case. 
Our task here is to figure out the conditions for the existence of a projective
logarithmic heat operator. 

Let $f:X\to S$ be a flat family of normal crossing varieties of dimension $d$
satisfying the assumpations of Proposition 1.3. Since $X$ is smooth, the (1.2)
will imply that $\lambda_{\alpha_1,...,\alpha_r}=0.$ Thus, for any local section
$D\in\SD^k_{X^{\dagger}/S^{\dagger}}$, we have
$$D=\sum_{\beta_1+\cdots+\beta_r\le k}
\lambda_{\beta_1,...,\beta_r}\p_1^{\beta_1}\cdots\p_r^{\beta_r}.$$
Namely, in this case, we have ($T_{X/S}(logY)$ isomorphic to the dual of  
$\Omega^1_{X/S}(logY)$) 
$$Sym^i_{\SO_X}(T_{X/S}(logY))\cong Gr_i(\SD_{X^{\dagger}/S^{\dagger}}),$$ 
which means that we have the canonical exact sequences (Proposition 1.4)
$$\CD
0 @>>>\SD^{k-1}_{X^{\dagger}} @>>>\SD^k_{X^{\dagger}} @>>>S^kT_X(logY) @>>>0 \\
@.   @AAA   @AAA   @AAA   @.      \\
0 @>>>\SD^{k-1}_{X^{\dagger}/S^{\dagger}}@>>>\SD^k_{X^{\dagger}/S^{\dagger}}
@>>>S^kT_{X/S}(logY)@>>> 0  \\
@.  @AAA  @AAA    @AAA  @.\\ 
@. 0 @. 0  @. 0  @.
\endCD\tag2.1$$

For convenience of applications, we summarize the discussions in section 1
for the special logarithmic varieties we consider in this paper. Let $f:Z\to T$
denote the logarithmic schemes: (1) $Z=X$ and $T=S$ with log structures
$logY$ and $log\Delta$, (2) $Z=X_s$ and $T=\{s\}$ with log structures $logY|_{X_s}$ 
and $log\Delta|_{\{s\}}$. Let $\Lambda$ and $\Lambda_i$ denote
the corresponding sheaf of logarithmic differential operators on $Z$. Then,
according to [Si], the following proposition (which is the definition of [Si])
assures that $\Lambda$ is 
$a$ $sheaf$ $of$ $split$ $almost$ $polynomial$ $rings$ $of$ $differential$
$operators$ $on$ $Z/T$. Thus $\Lambda$ and $\Lambda_i$ are endowed with
all the nice properties such as compatible with base changes and $\Lambda$ 
generated by $\Lambda_1$ as a ring (see [Si]). 

\proclaim{Proposition 2.1} The sheaf $\Lambda$ of rings of logarithmic 
differential operators on $Z$ over $T$ is a sheaf of $\SO_Z$-algebra
$\Lambda$ over $T$ with a filtration
$$\Lambda_0\subset\Lambda_1\subset\cdots\subset\Lambda_i\subset\cdots,$$
which satisfies\roster
\item $\Lambda=\bigcup^{\infty}{i=0}\Lambda_i$ and $\Lambda_i\cdot\Lambda_j
\subset\Lambda_{i+j}.$
\item The image of the morphism $\SO_Z\to\Lambda$ is equal to $\Lambda_0.$
\item The image of $f^{-1}(\SO_T)$ in $\SO_Z$ is contained in the center of
$\Lambda.$
\item The left and right $\SO_Z$-module structures on 
$Gr_i(\Lambda):=\Lambda_i/\Lambda_{i-1}$ are equal.
\item The sheaves of $\SO_Z$-modules $Gr_i(\Lambda)$ are coherent.
\item The sheaf of graded $\SO_Z$-algebra $Gr(\Lambda):=
\bigoplus^{\infty}_{i=0}Gr_i(\Lambda)$ is generated by $Gr_1(\Lambda)$ in the sense 
that the morphism of sheaves
$$Gr_1(\Lambda)\otimes_{\SO_Z}\cdots\otimes_{\SO_Z}Gr_1(\Lambda)\to Gr_i(\Lambda)$$ 
is surjective.
\item $\Lambda_0=\SO_Z$, $Gr_1(\Lambda)$ is locally free and $Gr(\Lambda)$ is
the symmetric algebra on $Gr_1(\Lambda)$.
\item There is a morphism $\xi:Gr_1(\Lambda)\to\Lambda_1$ of
left $\SO_Z$-modules splitting the projection $\Lambda_1\to Gr_1(\Lambda)$.
\endroster\endproclaim

\proclaim{Definition 2.1} For any coherent $\SO_X$-modules $E_1$ and $E_2$,
we define $$\SD^k_{X^{\dagger}/S^{\dagger}}(E_1,E_2):=
E_2\otimes_{\SO_X}\SD^k_{X^{\dagger}/S^{\dagger}}\otimes_{\SO_X}E_1^*$$
$$\SD^k_{X^{\dagger}}(E_1,E_2):=
E_2\otimes_{\SO_X}\SD^k_{X^{\dagger}}\otimes_{\SO_X}E_1^*,$$
where $E^*_1=\Cal Hom_{\SO_X}(E_1,\SO_X)$, and the notation
$\SD^k_{X^{\dagger}/S^{\dagger}}\otimes_{\SO_X}E_1^*$ (resp.
$\SD^k_{X^{\dagger}}\otimes_{\SO_X}E_1^*$ means that we use the right 
$\SO_X$-module structure of $\SD^k_{X^{\dagger}/S^{\dagger}}$ (resp.
$\SD^k_{X^{\dagger}}$). If $E_1=E_2=E$, we simply write 
$\SD^k_{X^{\dagger}/S^{\dagger}}(E)=
\SD^k_{X^{\dagger}/S^{\dagger}}(E_1,E_2)$ 
and $\SD^k_{X^{\dagger}}(E)=\SD^k_{X^{\dagger}}(E_1,E_2).$
\endproclaim

Let $\SL$ be a line bundle on $X$, and define the subsheaf $\SW_{X/S}(\SL)$
of $\SD^2_{X^{\dagger}}(\SL)$ to be
$$\SW_{X/S}(\SL):=\SD^1_{X^{\dagger}}(\SL)+
\SD^2_{X^{\dagger}/S^{\dagger}}(\SL).$$
Tensor firstly the sequence (2.1) on the right by $\SL^*$ as a right 
$\SO_X$-module then on the left by $\SL$ as a left $\SO_X$-module, we get
the commutative diagram for $k=2$
$$\CD
0 @>>>\SD^1_{X^{\dagger}}(\SL)@>>>\SD^2_{X^{\dagger}}(\SL)
@>\sigma>>S^2T_X(logY) @>>>0 \\
@.   @|   @AAA   @AAA   @.      \\
0 @>>>\SD^1_{X^{\dagger}}(\SL)@>>>\SW_{X/S}(\SL)
@>\sigma>>S^2T_{X/S}(logY)@>>> 0  \\
@.  @AAA  @AAA    @|  @.\\ 
0 @>>>\SD^1_{X^{\dagger}/S^{\dagger}}(\SL)@>>>
\SD^2_{X^{\dagger}/S^{\dagger}}(\SL)
@>\sigma>>S^2T_{X/S}(logY)@>>> 0,
\endCD$$
and the commutative diagram for $k=1$
$$\CD
@.   @.  @.                      f^*T_S(log\Delta) @.\\
@.    @.  @.                      @AAA              @. \\
0 @>>>\SO_X @>>>\SD^1_{X^{\dagger}}(\SL) @>\sigma_1>>T_X(logY) @>>>0 \\
@.   @|   @AAA   @AAA   @.      \\
0 @>>>\SO_X @>>>\SD^1_{X^{\dagger}/S^{\dagger}}(\SL)
@>\sigma_1>>T_{X/S}(logY)@>>> 0,  
\endCD\tag2.2$$
where the third vertical is the canonical exact sequence
$$0\to T_{X/S}(logY)\to T_X(logY)\to f^*T_S(log\Delta)\to 0.\tag2.3$$
Let $\varepsilon: \SD^1_{X^{\dagger}}(\SL)@>\sigma_1>> T_X(logY)\to f^*T_S(log\Delta)$ 
be the composition of canonical maps, one can see easily
from the diagram (2.2) that
$$ker(\varepsilon)=\SD^1_{X^{\dagger}/S^{\dagger}}(\SL).$$
Thus we have a surjection $\bar\varepsilon:\SW_{X/S}(\SL)\to f^*T_S(log\Delta)$
such that the following diagram is commutative
$$\CD
  f^*T_S(log\Delta)@= f^*T_S(log\Delta) @.  \\
   @A\varepsilon AA     @A\bar\varepsilon AA   @.      \\
\SD^1_{X^{\dagger}}(\SL)@>i>>\SW_{X/S}(\SL)
@>\sigma>>S^2T_{X/S}(logY)@>>> 0  \\
  @AAA  @AAA    @|  @.\\ 
\SD^1_{X^{\dagger}/S^{\dagger}}(\SL)@>>>
\SD^2_{X^{\dagger}/S^{\dagger}}(\SL)
@>\sigma>>S^2T_{X/S}(logY)@>>> 0
\endCD\tag2.4$$
Let $\SW_{X/S}(\SL)@>\sigma\oplus\bar\varepsilon>>S^2T_{X/S}(logY)\oplus
f^*T_S(log\Delta)$ be the surjection defined by
$\sigma\oplus\bar\varepsilon(D):=\sigma(D)\oplus\bar\varepsilon(D)$ for any
local section $D\in\SW_{X/S}(\SL)$. Then we have the exact sequence
$$0\to\SD^1_{X^{\dagger}/S^{\dagger}}(\SL)\to\SW_{X/S}(\SL)
@>\sigma\oplus\bar\varepsilon>>
S^2T_{X/S}(logY)\oplus f^*T_S(log\Delta)\to 0.\tag2.5$$

\proclaim{Definition 2.2} A logarithmic heat operator on $\SL$ over $S$ is an 
$\SO_S$-module homomorphism 
$$H: T_S(log\Delta)\to f_*\SW_{X/S}(\SL)\subset\SD^2_{X^{\dagger}}(\SL)$$
such that 
$$T_S(log\Delta)@>H>>f_*\SW_{X/S}(\SL)@>f_*\bar\varepsilon>>T_S(log\Delta)
\tag2.6$$
is the identity map. A logarithmic heat operator $H$ on $\SL$ is called
flat if $$H([\theta_1,\theta_2])=[H(\theta_1),H(\theta_2)]$$ for any local 
sections $\theta_1,\theta_2\in T_S(log\Delta)(U).$\endproclaim

Any $\SO_S$-linear map $\wt H: T_S(log\Delta)\to f_*\SW_{X/S}(\SL)/\SO_S$
has local lifting. Namely, there exists an open covering $\bigcup U=S$ such
that for each open set $U$ there is an $\SO_U$-linear map
$H_U:T_S(log\Delta)\to f_*\SW_{X/S}(\SL)|_U$ which reduces to $\wt H|_U$.

\proclaim{Definition 2.3} A projective logarithmic heat operator $\wt H$ 
on $\SL$ over $S$ is an $\SO_S$-linear map
$$\wt H: T_S(log\Delta)\to \frac{f_*\SW_{X/S}(\SL)}{\SO_S}$$
such that any local lifting $H_U$ is a logarithmic heat operator on
$\SL |_{f^{-1}(U)}$ over $U$. $\wt H$ is called projectively flat if any of the
local lifts $H_U$ satisfies 
$$H_U([\theta_1,\theta_2])=h_{\theta_1,\theta_2}+[H_U(\theta_1),H_U(\theta_2)]$$   
for some function $h_{\theta_1,\theta_2}\in\SO_S(V)$, where $V\subset U$ is any
open set of $U$ and $\theta_1,\theta_2\in T_S(log\Delta)(V)$.\endproclaim

In the following, we will figure out the conditions under which a projective
logarithmic heat operator on $\SL$ over $S$ do exist. As the same as in [GJ],
one can see that a (projective) logarithmic heat operator of $\SL$ over $S$ 
gives a (projective) logarithmic connection on $f_*\SL$. Firstly, it is clear 
that the map $$f_*\sigma:f_*\SW_{X/S}(\SL)\to f_*S^2T_{X/S}(logY)$$ factors 
through $f_*\SW_{X/S}(\SL)/\SO_S$, thus we have the map
$$\rho_{\wt H}:T_S(log\Delta)@>\wt H>>f_*\SW_{X/S}(\SL)/\SO_S@>f_*\sigma>>
f_*S^2T_{X/S}(logY),$$
which is called the $symbol$ of $\wt H$. By taking the direct image of 
$$0 @>>>\SO_X @>>>\SD^1_{X^{\dagger}/S^{\dagger}}(\SL)
@>\sigma_1>>T_{X/S}(logY)@>>> 0,\tag2.7$$
we have the connecting map $f_*T_{X/S}(logY)@>\cup[\SL]>>R^1f_*\SO_X$ and the map 
$$R^1f_*\SD^1_{X^{\dagger}/S^{\dagger}}(\SL)@>R^1f_*\sigma_1>>
R^1f_*T_{X/S}(logY)$$ 
induced by the symbol map 
$\SD^1_{X^{\dagger}/S^{\dagger}}(\SL)@>\sigma_1>>T_{X/S}(logY).$ Similarly,
from
$$0 @>>>\SD^1_{X^{\dagger}/S^{\dagger}}(\SL)@>>>
\SD^2_{X^{\dagger}/S^{\dagger}}(\SL)
@>\sigma>>S^2T_{X/S}(logY)@>>> 0,\tag2.8$$
we have the connecting map
$f_*S^2T_{X/S}(logY)@>c>>R^1f_*\SD^1_{X^{\dagger}/S^{\dagger}}(\SL)$ and thus
$$\mu_{\SL}:f_*S^2T_{X/S}(logY)@>c>>R^1f_*\SD^1_{X^{\dagger}/S^{\dagger}}(\SL)
@>R^1f_*\sigma_1>>
R^1f_*T_{X/S}(logY).$$ 
From the canonical exact sequence (2.3), we get the connecting map
$$\kappa_{X/S}:T_S(log\Delta)\to R^1f_*T_{X/S}(logY),$$ 
which is the Kodaira-Spencer map of the family $X/S$.

\proclaim{Theorem 2.1} Let $f:X\to S$, $\Delta$ and $Y=f^{-1}(\Delta)$ satisfy 
the assumptions of Proprosition 1.3 and $f_*\SO_X=\SO_S$, $\SL$ a line bundle on $X$. 
If there exists a symbol
$$\rho:T_S(log\Delta)\to f_*S^2T_{X/S}(logY)$$
such that the following two conditions hold\roster
\item $\mu_{\SL}\cdot\rho+\kappa_{X/S}=0$,
\item $f_*T_{X/S}(logY)@>\cup[\SL]>>R^1f_*\SO_X$ is an isomorphism.\endroster
Then there exists a unique projective logarithmic heat operator
$$\wt H:T_S(log\Delta)\to f_*\SW_{X/S}(\SL)/\SO_S$$
such that $\rho_{\wt H}=\rho$. In particular, there exists a projective
logarithmic connection on $f_*\SL$.\endproclaim

\demo{Proof} It is enough to prove that for any $\theta\in T_S(log\Delta)(U)$
there exists a unique lifting of $\rho(\theta)\oplus\theta$ to
$f_*\SW_{X/S}(\SL)(U)$ up to a section of $\SO_S(U)$. Thus we consider the
commutative diagram
$$\CD
T_{X/S}(logY)@>>>\frac{\SW_{X/S}(\SL)}{\SO_X}
@>\sigma\oplus\bar\varepsilon>>S^2T_{X/S}(logY)\oplus f^*T_S(log\Delta)\\
 @A\sigma_1AA  @AAA   @|      \\
\SD^1_{X^{\dagger}/S^{\dagger}}(\SL)@>>>\SW_{X/S}(\SL)
@>\sigma\oplus\bar\varepsilon>>
S^2T_{X/S}(logY)\oplus f^*T_S(log\Delta)\\
@AAA       @AAA  @. \\
\SO_X @= \SO_X, @. \endCD$$
which gives the induced commutative diagram
$$\CD 
R^1f_*\SO_X@=R^1f_*\SO_X @. \\
@A\cup[\SL]AA    @A\nu AA   @.\\
f_*T_{X/S}(logY)@>>>f_*\frac{\SW_{X/S}(\SL)}{\SO_X}
@>>>f_*S^2T_{X/S}(logY)\oplus T_S(log\Delta)@>o>> \\
 @AAA  @AAA   @|   @.   \\
f_*\SD^1_{X^{\dagger}/S^{\dagger}}(\SL)@>>>f_*\SW_{X/S}(\SL)
@>>> f_*S^2T_{X/S}(logY)\oplus T_S(log\Delta)@>>> \\
  @AAA       @AAA  @. \\
 f_*\SO_X @= f_*\SO_X, @. \endCD$$
where $f_*S^2T_{X/S}(logY)\oplus T_S(log\Delta)@>o>>R^1f_*T_{X/S}(logY)$
is the connecting map. We claim that
$$o(\rho(\theta)\oplus\theta)=\mu_{\SL}\cdot\rho(\theta)+\kappa_{X/S}(\theta).$$
If it is true, by the condition (1), we will have a lifting 
$\wt H_U(\theta)\in f_*\frac{\SW_{X/S}(\SL)}{\SO_X}(U).$ 
By the surjectivity in condition (2), there exists a section 
$s\in f_*T_{X/S}(logY)(U)$ such that 
$s\cup [\SL]=\nu(\wt H_U(\theta))$.
Thus there exists a $H_U(\theta)\in f_*\SW_{X/S}(\SL)(U)$ such that
$\wt{H_U(\theta)}=\wt H_U(\theta)-s,$ which is also a lifting of 
$\rho(\theta)\oplus\theta$. The injectivity in condition (2) implies that
such $H_U(\theta)$ is unique up to a section of $f_*\SO_X(U)=\SO_S(U)$. This
gives a unique projective logarithmic heat operator
$$\wt H: T_S(log\Delta)\to\frac{f_*\SW_{X/S}(\SL)}{\SO_S}.$$

Now we show the claim by considering the following commutative diagrams 
$$\CD
0@>>>\frac{\SD^1_{X^{\dagger}/S^{\dagger}}(\SL)}{\SO_X}@>>>
\frac{\SW_{X/S}(\SL)}{\SO_X}@>\sigma\oplus\bar\varepsilon>>
S^2T_{X/S}(logY)\oplus f^*T_S(log\Delta)\\
@.  @| @AAA    @AAA   \\ 
0@>>>\frac{\SD^1_{X^{\dagger}/S^{\dagger}}(\SL)}{\SO_X}@>>>
\frac{\SD^2_{X^{\dagger}/S^{\dagger}}(\SL)}{\SO_X}@>\sigma>>
S^2T_{X/S}(logY)
\endCD$$

$$\CD
0@>>>\frac{\SD^1_{X^{\dagger}/S^{\dagger}}(\SL)}{\SO_X}@>>>
\frac{\SW_{X/S}(\SL)}{\SO_X}@>\sigma\oplus\bar\varepsilon>>
S^2T_{X/S}(logY)\oplus f^*T_S(log\Delta)\\
@.  @| @A i AA    @AAA   \\ 
0@>>>\frac{\SD^1_{X^{\dagger}/S^{\dagger}}(\SL)}{\SO_X}@>>>
\frac{\SD^1_{X^{\dagger}}(\SL)}{\SO_X}@>\varepsilon>>
f^*T_S(log\Delta)\\
@. @V\sigma_1VV  @V\sigma_1VV  @|  \\
0@>>>T_{X/S}(logY)@>>>T_X(logY)@>>>f^*T_S(log\Delta),
\endCD$$
from which we have commutative diagrams for the connecting maps
$$\CD
f_*S^2T_{X/S}(logY)\oplus T_S(log\Delta)@>o>>
R^1f_*\frac{\SD^1_{X^{\dagger}/S^{\dagger}}(\SL)}{\SO_X}
\cong R^1f_*T_{X/S}(logY)\\
@AAA  @| @. \\
f_*S^2T_{X/S}(logY)@>\mu_{\SL}>>
R^1f_*\frac{\SD^1_{X^{\dagger}/S^{\dagger}}(\SL)}{\SO_X}
\cong R^1f_*T_{X/S}(logY)
\endCD$$
and
$$\CD
f_*S^2T_{X/S}(logY)\oplus T_S(log\Delta)@>o>>R^1f_*
\frac{\SD^1_{X^{\dagger}/S^{\dagger}}(\SL)}{\SO_X}\cong R^1f_*T_{X/S}(logY)\\
@AAA    @| @. \\
T_S(log\Delta)@>\kappa_{X/S}>>
R^1f_*\frac{\SD^1_{X^{\dagger}/S^{\dagger}}(\SL)}{\SO_X}\cong R^1f_*T_{X/S}(logY).
\endCD$$
Thus the claim  $o(\rho(\theta)\oplus\theta)=
\mu_{\SL}\cdot\rho(\theta)+\kappa_{X/S}(\theta)$ is indeed true. 
\enddemo

\remark{Remark 2.1} From the proof, we see that the local lifting exists
when the map in condition (2) is surjective, and the injectivity was only
used to assure the uniqueness of the local lifting. Thus in some cases the
map in condition (2) is only surjective but one has a natural way to choose
the lifting uniquely, we still have the heat operator. For example, the map
in [GJ] is zero but one can choose uniquely the $\Cal G$-invariant lifting.
\endremark

We are now going to describe the maps $\cup [\SL]$ and $\mu_{\SL}$. For any
$s\in\Delta$, $\{s\}=\text{Spec}k(s)$ has the induced log structure, we denote
this logarithmic point by $s^{\dagger}$, then the logarithmic fibre 
$X_s^{\dagger}$ over $s^{\dagger}$ is $(X_s,logY|_{X_s}).$ Thus
$$\Omega^1_{X_s}(log)=\Omega^1_{X/S}(logY)|_{X_s},\quad 
\SD^k_{X_s^{\dagger}}=\SD^k_{X^{\dagger}/S^{\dagger}}|_{X_s}$$
and, if the dimensions of $H^0(T_{X_s}(log))$, $H^0(S^2T_{X_s}(log))$,
$H^1(\SO_{X_s})$, $H^1(T_{X_s}(log))$ are constant (for $s$), then fibrewisely
the maps $\cup [\SL]$ and $\mu_{\SL}$ are the following maps
$$\cup[\SL_s]:H^0(T_{X_s}(log))\to H^1(\SO_{X_s})$$
$$\mu_{\SL_s}:H^0(S^2T_{X_s}(log))\to H^1(T_{X_s}(log))$$
where $\SL_s=\SL|_{X_s}$ and $\cup[\SL_s]$ is the connecting map of
$$0\to\SO_{X_s}\to\SD^1_{X_s^{\dagger}}(\SL_s)@>\sigma_1>>T_{X_s}(log)\to 0\tag2.9$$
and $\mu_{\SL_s}$ is the connecting map 
$H^0(S^2T_{X_s}(log))\to H^1(\SD^1_{X_s^{\dagger}}(\SL_s))$ of
$$0\to\SD^1_{X_s^{\dagger}}(\SL_s)\to\SD^2_{X_s^{\dagger}}(\SL_s)@>\sigma_2>>
S^2T_{X_s}(log)\to 0,\tag2.10$$
composing with the natural map 
$H^1(\SD^1_{X_s^{\dagger}}(\SL_s))@>H^1(\sigma_1)>>H^1(T_{X_s}(log)).$

Let $[\SL_s]\in H^1(\Omega^1_{X_s}(log))$ denote the extension class of (2.9),
then the map $\cup[\SL_s]$ means the cup product. In general, for any class
$cl\in H^1(\Omega^1_{X_s}(log))$, one has the natural cup product map
$$H^0(\otimes^kT_{X_s}(log))@>\cup cl>>H^1(\otimes^{k-1}T_{X_s}(log))$$
and for any $\omega\in H^0(S^kT_{X_s}(log))$ the symbol $\omega\cup cl$ means
that we consider $\omega$ as a symmetric tensor. For any line bundle $L$ on
$X_s$, we define the Chern class $c_1(L)\in H^1(\Omega^1_{X_s}(log))$ of $L$ to 
be the image of usual
Chern class of $L$ under the natural map
$H^1(\Omega^1_{X_s})\to H^1(\Omega^1_{X_s}(log)).$
More precisely, let $\bar d:\SO_{X_s}@>d>>\Omega^1_{X_s}\to\Omega^1_{X_s}(log)$
and $dl:\SO_{X_s}^*\to\Omega^1_{X_s}(log)$ be defined as $dl(u)=\frac{1}{u}\bar du$. 
Then $dl$ is a morphism of abelian sheaves and induces a morphism 
$$H^1(\SO_{X_s}^*)@>c_1>>H^1(\Omega^1_{X_s}(log))$$
of abelian groups. The Chern class $c_1(L)$ of $L$ is defined to be the image 
of this morphism. With these notation, we have

\proclaim{Proposition 2.2} 
The extension class $[\SL_s]\in H^1(\Omega^1_{X_s}(log))$ of
$$0\to\SO_{X_s}\to\SD^1_{X_s^{\dagger}}(\SL_s)@>\sigma_1>>T_{X_s}(log)\to 0\tag2.11$$
is equal to the Chern class $c_1(\SL_s)$ and for any $\omega\in H^0(S^2T_{X_s}(log))$, 
we have
$$\mu_{\SL_s}(\omega)=-\omega\cup c_1(\SL_s)+\mu_{\SO_{X_s}}(\omega).$$
\endproclaim

\demo{Proof} We check firstly the following descriptions about the symbol
maps.\roster\item If $D$ is a local section of $\SD^1_{X_s^{\dagger}}(\SL_s)$,
its image $\sigma_1(D)\in T_{X_s}(log)$ is determined by the requirement,
for all $a\in\SO_{X_s}$ and all $s\in\SL_s$
$$<\sigma_1(D),\bar da>\cdot s=\sigma_1(D)(a)\cdot s=D(a\cdot s)-a\cdot D(s).
\tag2.12$$
\item If $D$ is a local section of $\SD^2_{X_s^{\dagger}}(\SL_s)$, its image
$\sigma_2(D)\in S^2T_{X_s}(log)$ is characterized by the formula, for all
$a,\,b\in\SO_{X_s}$ and all $s\in\SL_s$
$$\sigma_2(D)(a,b)\cdot s=D(ab\cdot s)-a\cdot D(b\cdot s)-b\cdot D(a\cdot s)
+ab\cdot D(s).\tag2.13$$\endroster
(1) is clear by the definition $\sigma_1(D)(a)\cdot s=[D,a](1)\cdot s
=D(a\cdot s)-a\cdot D(s)$ in Proposition 1.2 (1). To check (2), one can write
$D=D_{<2}+\sum\lambda_{ij}\p_i\p_j$, where $D_{<2}$ denotes the part with order
smaller than $2$. Then, by definition in Proposition 1.4,
$$<\sigma_2(D),\bar da\otimes\bar db>\cdot s:=\sigma_2(D)(a,b)=
(\sum\lambda_{ij}[\p_i,a][\p_j,b]+ \sum\lambda_{ij}[\p_i,b][\p_j,a])(1)\cdot s.$$
Thus (2) is clear from the following computations 
$$\aligned\sum\lambda_{ij}[\p_i,a][\p_j,b]&=\sum\lambda_{ij}\p_ia 
\p_jb-\sum\lambda_{ij}a\p_i\p_jb-\sum\lambda_{ij}\p_iab\p_j+\sum\lambda_{ij}a
\p_ib\p_j\\&=Dab-aDb-bDa+abD+b[D_{<2},a]-[D_{<2},a]b\\&
-[\sum\lambda_{ij}\p_i[\p_j,a],b]-[\sum\lambda_{ij}[\p_i,b]\p_j,a]+
\sum\lambda_{ij}[\p_i,b][\p_j,a]\\&=Dab-aDb-bDa+abD
-\sum\lambda_{ij}[\p_i,b][\p_j,a].\endaligned$$

Let $\SU=\{U_i\}_{i\in I}$ be an affine open cover of $X_s$ trivializing
$\SL_s$ and $s_i:\SO_{U_i}\cong\SL_s|_{U_i}$, $s_j=u_{ij}\cdot s_i$ on
$U_{ij}=U_i\cap U_j$. Then $c_1(\SL_s)\in H^1(\Omega^1_{X_s}(log))$ is given 
by the 1-cocycle 
$$\{\frac{\bar du_{ij}}{u_{ij}}\}\in C^1(\SU,\Omega^1_{X_s}(log)).$$
The sequence (2.11) is locally splitting, and there exist morphisms of
$\SO_{U_i}$-modoules
$$\rho_i:T_{X_s}(log)(U_i)\to\SD^1_{X_s^{\dagger}}(\SL_s)(U_i)$$
such that $\sigma_1\circ\rho_i(\theta)=\theta$ for any $\theta\in T_{X_s}(log)(U_i).$ 
Let $\rho_i(\theta)(s_i)=\omega_i(\theta)\cdot s_i$, then
$\omega_i\in T_{X_s}(log)(U_i)^*=\Omega^1_{X_s}(log)(U_i)$ since $\rho_i$
is a morphism of $\SO_{U_i}$-modules. For any $\theta\in T_{X_s}(log)(U_{ij})$
$$\rho_i(\theta)(s_i)-\rho_j(\theta)(s_i)=<\theta,\omega_i-\omega_j+
\frac{\bar du_{ij}}{u_{ij}}>\cdot s_i.$$
Thus the extension class of (2.11) is $c_1(\SL_s)$.

Given $\omega\in H^0(S^2T_{X_s}(log))$, let 
$D_i\in H^0(U_i,\SD^2_{X_s^{\dagger}}(\SO_{X_s}))$ be the lifting of
$\omega_i=\omega|_{U_i}$, thus 
$$\{D_{ij}\}_{i<j}=\{D_j-D_i\}_{i<j}\in C^1(\SU,\SD^1_{X_s^{\dagger}}(\SO_{X_s})).$$
Then $\mu_{\SO_{X_s}}(\omega)\in H^1(T_{X_s}(log))$ is given by the $1$-cocycle
$\{v_{ij}\}\in C^1(\SU, T_{X_s}(log))$, where $v_{ij}\in H^0(U_{ij},T_{X_s}(log))$
is the image $\sigma_1(D_{ij})$ of $D_{ij}=D_j-D_i$, namely, for all
$a\in \SO_{X_s}(U_{ij})$, 
$$<v_{ij},\bar da>=v_{ij}(a)=D_{ij}(a)-aD_{ij}(1).$$

To compute $\mu_{\SL_s}(\omega)$, we see that by definition 
$\wt D_i=s_i\otimes D_i\otimes  s_i^*\in H^0(U_i,\SD^2_{X_s^{\dagger}}(\SL_s)$ 
is a lifting of
$\omega_i$, and thus $\wt D_i(a\cdot s_i)=D_i(a)\cdot s_i$ and
$$\{\wt D_{ij}\}_{i<j}=
\{\wt D_j-\wt D_i\}_{i<j}\in C^1(\SU,\SD^1_{X_s^{\dagger}}(\SL_s)).$$
Then $\mu_{\SL_s}(\omega)\in H^1(T_{X_s}(log))$ is given by the $1$-cocycle
$\{\wt v_{ij}\}\in C^1(\SU, T_{X_s}(log))$, 
where $\wt v_{ij}\in H^0(U_{ij},T_{X_s}(log))$
is the image $\sigma_1(\wt D_{ij})$ of $\wt D_{ij}$, namely, for all
$a\in \SO_{X_s}(U_{ij})$, by using (2.13), we have 
$$\aligned<\wt v_{ij},\bar da>\cdot s_j&=\wt v_{ij}(a)=
\wt D_{ij}(a\cdot s_j)-a\cdot\wt D_{ij}(s_j)\\
&=D_j(a)\cdot s_j-D_i(au_{ij})\cdot s_i-aD_j(1)\cdot s_j+aD_i(U_{ij})\cdot s_i\\
&=(v_{ij}(a)-<\omega,\bar da\otimes\frac{\bar du_{ij}}{u_{ij}}>)\cdot s_j.
\endaligned$$
Hence $\wt v_{ij}=v_{ij}-\omega\cup\frac{\bar du_{ij}}{u_{ij}}$ and
$\mu_{\SL_s}(\omega)=\mu_{\SO_{X_s}}(\omega)-\omega\cup c_1(\SL_s).$

\enddemo

\heading \S 3 Logarithmic operators on generalized Jacobians\endheading

In this section, we will verify the conditions in Theorem 2.1 for a family
of generalized Jacobians of stable curves, and thus show the existence of
logarithmic heat operator. Let $(\Cal C,\Cal C_{\Delta})\to (S,\Delta)$ be 
a flat family of stable
curves satisfying the assumpations in Proposition 1.3, namely, $S$, $\Cal C$
are regular schemes and $\Delta$ a (reduced) normal crossing divisor. It is
well known there exists a projective $S$-scheme $f: J(\Cal C)\to S$ such that
for any $s\in S$ the fibre $J(\Cal C)_s$ is the generalized Jacobian 
$J(\Cal C_s)$ of $\Cal C_s$.

\proclaim{Lemma 3.1} If  $(\Cal C,\Cal C_{\Delta})\to (S,\Delta)$  satisfies
the assumpations of Proposition 1.3, then so do $f:J(\Cal C)\to S$.\endproclaim

\demo{Proof} By deformation theory of torsion free sheaves with rank one, for any point
$y\in J(\Cal C)$ corresponds to a torsion free sheaf $\SF$ on $\Cal C_{f(y)}$
such that $\SF$ is not locally free at a double point $x\in\Cal C_{f(y)}$,
there are integers $l_1$, $l_2$ such that
$$\hat\SO_{J(\Cal C),y}[[u_1,...,u_{l_1}]]\cong\hat\SO_{\Cal C,x}[[v_1,...,v_{l_2}]].$$
Thus $f:J(\Cal C)\to S$ satisfies the assumpations in Proposition 1.3 if
$\Cal C/S$ satisfies them. In particular, $J(\Cal C)$ is regular and all
fibres $J(\Cal C_s)$ are normal crossing varieties.\enddemo 

For simplicity, we assume that all fibres $\Cal C_s=C$ ($s\in\Delta$) are 
irreducible and smooth except one node $x_0$. We recall briefly some facts 
about the so called generalized Jacobian $J^d(C)$ (we write $J(C)$ for $J^0(C)$)
of a projective (singular) curve $C$ of (arithmetic) genus $g$.

The $J^d(C)$ is defined to be the 
moduli space of rank one torsion free sheaves with degree $d$. There is a
natural ample line bundle 
$$\Theta_{J^d(C)}=detH^*(\SF)^{-1}\otimes (det\SF_y)^{d+1-g}$$ 
called theta line bundle on $J^d(C)$, where $\SF$ is a universal family on 
$C\times J^d(C)$, $detH^*(\SF)$ is the
determinant of cohomology and $\SF_y$ denotes the restriction of $\SF$ to
$\{y\}\times J^d(C)$ for a fixed smooth point $y\in C$.
This construction can be generalized to relative case, namely, for a family
of curves $\Cal C/T$, one can construct a family of generalized Jacobians
$J^d(\Cal C)/T$ and a line bundle $\Theta$ on $J^d(\Cal C)$ such that
each fibre $J^d(\Cal C)_t$ is the generalized Jacobian $J^d(\Cal C_t)$ and the
restriction of $\Theta$ to $J^d(\Cal C)_t$ is the theta line bundle 
$\Theta_{J^d(\Cal C_t)}$.

Let $\pi:\wt C\to C$ be the normalization and $\pi^{-1}(x_0)=\{x_1,x_2\}$, let
$P=\Bbb P(\SE_{x_1}\oplus\SE_{x_2})$ and 
$\SE_{x_1}\oplus\SE_{x_2}\to\SO(1)\to 0$
be the universal quotient on $P$, where $\SE$ is a univeral line bundle over
$\wt C\times J^d(\wt C)$. We consider the diagram
$$\CD
P=\Bbb P(\SE_{x_1}\oplus\SE_{x_2})@>\phi>>J^d(C)\\
@V\rho VV             @. \\
J^d(\wt C)       @.             \endCD$$
where $\rho$ is the natural projection, and $\phi$ is defined as follows:
for any $(L,q):=(L,L_{x_1}\oplus L_{x_2}@>q>>\Bbb C)\in P$, $\phi(L,q)$ is the
kernel of $\pi_*L@>q>>\,_{x_0}\Bbb C\to 0.$ 

\proclaim{Lemma 3.2} Let $W\subset J^d(C)$ be the reduced subscheme of 
non-locally free sheaves and $D_1$, $D_2$ be the sections of 
$P@>\rho>>J^d(\wt C)$ given by projections 
$\SE_{x_1}\oplus\SE_{x_2}\to\SE_{x_1}$ and $\SE_{x_1}\oplus\SE_{x_2}\to\SE_{x_2}$. Then
\roster
\item $P@>\phi>>J^d(C)$ is the normalization of $J^d(C)$, and $W$ is the
non-normal locus of $J^d(C)$.
\item $\phi^{-1}(W)=D_1+D_2$ and $\phi|_{D_i}:D_i\to W$ ($i=1,2$) are isomorphisms.
\item For any integer $k>0$, $\Theta_P:=\phi^*(\Theta^k_{J^d(C)})=
\SO(1)^k\otimes\rho^*\SE_y^{-k}\otimes\rho^*\Theta^k_{J^d(\wt C)}$ and
$$\rho_*\SO(1)^k=\bigoplus_{j=0}^k\SE_{x_1}^j\otimes\SE_{x_2}^{k-j}.$$
\endroster\endproclaim

\demo{Proof} This is the special case (rank one) of [NR] and [Su].\enddemo

\proclaim{Lemma 3.3} Fix a line bundle $\SL=\SO_{\Bbb P^1}(1)$ and 
two points $p_1=(1,0)$, $p_2=(0,1)$ of $X=\Bbb P^1$,
which give a logarithmic structure on $X$. For any $D\in H^0(\SD^1_{X^{\dagger}}(\SL))$, 
if there exist $c\in\Bbb C^*$ and
a nonzero section $s\in H^0(\SL)$ satisfying $s(p_2)=c\cdot s(p_1)$, such that
$$D(s)(p_2)=c\cdot D(s)(p_1).$$
Then the symbol of $D$ is trivial.\endproclaim

\demo{Proof} $\Bbb P^1$ is covered by $V_1=\Bbb P^1\ssm\{p_2\}
=Spec\Bbb C[\frac{x_2}{x_1}]$ and $V_2=\Bbb P^1\ssm\{p_1\}=Spec\Bbb C[\frac{x_1}{x_2}]$,
and there is a global vector field $\p\in H^0(T_X(-p_1-p_2))$ such that
$$\p_1:=\p|_{V_1}=u\frac{\p\,}{\p u},\quad\p_2:=\p|_{V_2}=\,-v\frac{\p\,}{\p v}$$ 
where $u=\frac{x_2}{x_1}$, $v=\frac{x_1}{x_2}=\frac{1}{u}$. The space
$H^0(T_X(-p_1-p_2))$ is generated by $\p$.

We see that $\SL|_{V_1}=\Bbb C[u]\cdot x_1$ and $\SL|_{V_2}=\Bbb C[v]\cdot x_2$, 
thus any section $s\in H^0(\SL)$ has the form
$$s_1:=s|_{V_1}=(a_0+a_1u)\cdot x_1,\quad s_2:=s|_{V_2}=(a_1+a_0v)\cdot x_2,$$
where $(a_0,a_1)=(s(p_1),s(p_2))\in\Bbb C^2$. Therefore, for any
$D\in H^0(\SD^1_{X^{\dagger}}(\SL))$, there exists $(b_0,b_1)\in\Bbb C^2$
such that $$D(s)|_{V_1}=(b_0+b_1u)\cdot x_1,\quad D(s)|_{V_2}=(b_1+b_0v)\cdot x_2.$$
If the symbol of $D$ is $k\cdot\p$, by the definition of symbol, we have
$$\aligned &D(s)|_{V_1}=D(s_1)=(a_0+a_1u)\cdot D(x_1)+ka_1u\cdot x_1\\
 &D(s)|_{V_2}=D(s_2)=(a_1+a_0v)\cdot D(x_2)-ka_0v\cdot x_2 .\endaligned$$
Thus $D$ is determined by any given number $D(1)\in\Bbb C$ such that 
 $$D(x_1)=D(1)\cdot x_1,\quad D(x_2)=(D(1)+k)\cdot x_2,$$
and one checks that for any given number $D(1)\in\Bbb C$ the above definition
gives indeed a global differential operator of $\SL$ with symbol $k\cdot\p$. 
It is easy to see that for any $s\in H^0(\SL)$ and $c\in\Bbb C^*$ 
$$D(s)(p_2)-c\cdot D(s)(p_1)=(s(p_2)-c\cdot s(p_1))D(1)+k\cdot s(p_2).$$
Thus, if there exist a nonzero $s$ and $c$ such that $s(p_2)=c\cdot s(p_1)$,
we have 
$$D(s)(p_2)-c\cdot D(s)(p_1)=ks(p_2),$$
which is nonzero except $k=0$ since $s(p_2)\neq 0$ (otherwise $s(p_1)=0$ and
$s$ will have at least two zero points).\enddemo 

The fact that $X=J^d(C)$ is a degenerating fibre of flat family means that
$X$ is more special than usual normal crossing varieties. For example, its
cohomology has low bound ($h^1(\SO_X)\ge g$) and there is a logarithmic structure
on it. Moreover, we have

\proclaim{Proposition 3.1} Let $X=J^d(C)$ be the moduli space of torsion free
sheaves on $C$ with rank one and degree $d$, and $\SL$ be the theta line bundle
on $X$. Then for any logarithmic structure on $X$ in the sense of [KN] and any
integer $k>0$
$$H^0(\SD^1_{X^{\dagger}}(\SL^k))\cong\Bbb C,\quad H^0(T_X(log))\cong H^1(\SO_X)\cong
\Bbb C^g.$$\endproclaim

\demo{Proof} It is enough to prove Proposition 3.1 for $k=1$ since 
$c_1(\SL^k)=kc_1(\SL)$ and thus we have isomorphism 
$\SD^1_{X^{\dagger}}(\SL)\cong \SD^1_{X^{\dagger}}(\SL^k)$ by Proposition 2.2.

Local computation shows that $\phi^*\Omega^1_X(log)\cong\Omega^1_P(log(D_1+D_2))$ 
and thus the natural map 
$T_{P^{\dagger}}:=T_P(log(D_1+D_2))\to \phi^*T_X(log)$
is an isomorphism, where $P^{\dagger}=(P, log(D_1+D_2))$, and the diagram 
$$\CD
0@>>>\SO_P@>>>\SD^1_{P^{\dagger}}(\phi^*\SL)@>>> T_{P^{\dagger}}@>>> 0\\
@. @|  @VVV  @VVV @.\\
0@>>>\phi^*\SO_X@>>>\phi^*\SD^1_{X^{\dagger}}(\SL)@>>>\phi^*T_X(log)@>>>0
\endCD$$
implies that $\SD^1_{P^{\dagger}}(\Theta_P)\to\phi^*\SD^1_{X^{\dagger}}(\SL)$
is an isomorphism. Hence we have  
$$\CD
0@>>>\phi_*\SO_P@>>>\phi_*\SD^1_{P^{\dagger}}(\Theta_P)@>>>\phi_* T_{P^{\dagger}}@>>> 0\\
@. @AAA  @AAA  @AAA @.\\
0@>>>\SO_X@>>>\SD^1_{X^{\dagger}}(\SL)@>>>T_X(log)@>>>0
\endCD$$
Any operator $D\in H^0(\SD^1_{X^{\dagger}}(\SL))$ with nonzero symbol will
give an operator $D\in H^0(\SD^1_{P^{\dagger}}(\Theta_P))$ with nonzero symbol
$\sigma(D)\in H^0(T_{P^{\dagger}})$. On the other hand, we have
$$\CD
0@>>>\SO_P@>>>\SD^1_P(\rho^*\Theta_{J^d(\wt C)})@>>>T_P@>>> 0\\
@. @|  @VVV  @VVV @.\\
0@>>>\rho^*\SO_{J^d(\wt C)} @>>>\rho^*\SD^1_{J^d(\wt C)}(\Theta_{J^d(\wt C)})@>>>
\rho^*T_{J^d(\wt C)}@>>>0
\endCD$$ 
and $\SO(D_1)=\SO(1)\otimes\rho^*\SE_{x_1}^{-1}$. It is easy to see that
$$\SD^1_{P^{\dagger}}(\Theta_P)\hookrightarrow\SD^1_P(\Theta_P)\hookrightarrow
\SD^1_P(\Theta_P)\otimes\SO(D_1)=\SO(D_1)\otimes
\SD^1_P(\rho^*\wt\SL),$$
where $\wt\SL=\Theta_{J^d(\wt C)}\otimes\SE_y^{-1}\otimes\SE_{x_1}$, which is
isomorphic to $\Theta_{J^d(\wt C)}$. Consider 
$$0\to T_{P/J^d(\wt C)}(log)\to T_{P^{\dagger}}\to\rho^*T_{J^d(\wt C)},$$
if the image of $\sigma(D)$ in $H^0(\rho^*T_{J^d(\wt C)})$ is nonzero, then
the connecting map of
$$0@>>>\SO(D_1)\otimes\rho^*\SO_{J^d(\wt C)} @>>>\SO(D_1)
\otimes\rho^*\SD^1_{J^d(\wt C)}(\Theta_{J^d(\wt C)})
@>>>\SO(D_1)\otimes\rho^*T_{J^d(\wt C)}@>>>0$$
is not injective, which is impossible since the dimension of
$$\aligned &H^0(\SO(D_1)\otimes\rho^*\SD^1_{J^d(\wt C)}(\Theta_{J^d(\wt C)}))\\
&=H^0(\SD^1_{J^d(\wt C)}(\Theta_{J^d(\wt C)}))\oplus
H^0(\SD^1_{J^d(\wt C)}(\Theta_{J^d(\wt C)},\Theta_{J^d(\wt C)}
\otimes\SE_{x_2}\otimes\SE_{x_1}^{-1}))\endaligned$$
is two (see [We]) and the space 
$$H^0(\SO(D_1)\otimes\rho^*\SO_{J^d(\wt C)})=H^0(\SO_{J^d(\wt C)})\oplus
H^0(\SO_{J^d(\wt C)}\otimes\SE_{x_1}^{-1}\otimes\SE_{x_2})$$
is also of dimension two, where we remark that on a Jacobian 
$\Theta_{J^d(\wt C)}\otimes\SE_{x_2}\otimes\SE_{x_1}^{-1}$ and 
$\SO_{J^d(\wt C)}\otimes\SE_{x_1}^{-1}\otimes\SE_{x_2}$ are isomorphic to 
$\Theta_{J^d(\wt C)}$
and $\SO_{J^d(\wt C)}$. Thus if $D\in H^0(\SD^1_{X^{\dagger}}(\SL))$ has nonzero
symbol, then $D$ gives an operator 
$D\in H^0(\SD^1_{P^{\dagger}/J^d(\wt C)}(\Theta_P))$, which induces an operator 
$D\in H^0(\SD^1_{F^{\dagger}}(\wt L))$
with nonzero symbol for general fibre $F=\Bbb P^1$ of $\rho:P\to J^d(\wt C)$ and
the $\wt L=\SO_{\Bbb P^1}(1)$. By using Lemma 3.3, this is impossible. In fact, 
for any fibre $F$, let $p_1=F\cap D_1$, $p_2=F\cap D_2$ and $C_F=\phi(F)$, then
$\phi_F: F\to C_F$ is the map identifying $p_1$, $p_2$ into a node $p_0$ of $C_F$
and $X$ is covered by $\{C_F\}$. Let $L=\SL|_{C_F}$, then 
$\phi_F^*L=\phi^*(\SL)|_F=\SO_{\Bbb P^1}(1)$. Thus $L$ is defined by 
$$0\to L\to\phi_{F_*}\wt L@>c>>\,_{p_0}\Bbb C\to 0,$$
where $c\in\Bbb C^*$ and the quotient map $\wt L_{x_1}\oplus\wt L_{x_2}@>c>>\Bbb C$ 
is defined by $c(a,b)=b-c\cdot a$. Choose $F$ such that $H^0(\SL|_{C_F})\neq 0$, 
then we find a $c\in \Bbb C^*$ and a nonzero $s\in H^0(\SO_{\Bbb P^1}(1))$
satisfying that $s(p_2)=c\cdot s(p_1)$. Since $D$ has to induce a morphism 
$\SL\to\SL$ of abelian group sheaves, $D(s)$ has
to satisfy $D(s)(p_2)=c\cdot D(s)(p_1)$, which means that $D$ has zero symbol
by Lemma 3.3.

To see that $H^0(T_X(log))=H^1(\SO_X)=\Bbb C^g$, we remark that both spaces
have at least dimension $g$, then we only need to check that $dim\,H^1(\SO_X)\le g$. 
This is easy to see by using
$$0\to\SO_X\to\phi_*\SO_P\to\SO_W\to 0$$
and $H^1(\phi_*\SO_P)=H^1(\SO_P)=H^1(\rho_*\SO_P)=H^1(\SO_{J^d(\wt C)})=\Bbb C^{g-1}.$  
\enddemo

\proclaim{Lemma 3.4} For any logarithmic structure on $X=J^d(C)$ in the sense
of [KN], we have $H^0(S^2T_X(log))=S^2H^0(T_X(log))$.\endproclaim

\demo{Proof} It is enough to show that 
$$h^0(S^2T_X(log)):=dim\,H^0(S^2T_X(log))\le dim\,S^2H^0(T_X(log))=\frac{g(g+1)}{2}.$$
To prove it, let $\SF:=\phi^*T_X(log)=T_P(log(D_1+D_2))$, 
$\SF':=T_{P/J^d(\wt C)}(-D_1-D_2)$, $\SF":=\rho^*T_{J^d(\wt C)}$ and use the exact squence
$$0\to\SF'\to\SF\to\SF"\to 0,$$
one has $h^0(S^2T_X(log))\le h^0(S^2\SF)$ and the following two exact sequences
$$0\to\Cal G\to S^2\SF\to S^2(\SF")\to 0,$$
$$0\to S^2(\SF')\to\Cal G\to \SF'\otimes\SF"\to 0.$$
Thus, by using $h^0(S^2(\SF"))=h^0(S^2T_{J^d(\wt C)})=\frac{g(g-1)}{2}$, we have
$$h^0(S^2T_X(log))\le\frac{g(g-1)}{2}+h^0(S^2\SF')+h^0(T_{J^d(\wt C)}\otimes\rho_*\SF').$$
To compute $\SF'$, noting that $\SO_P(D_i)=
\SO(1)\otimes\rho^*\SE_{x_i}^{-1}$ and using the exact sequence
$$0\to\SO_P\to\SO(1)\otimes\rho^*(\SE_{x_1}^{-1}\oplus\SE_{x_2}^{-1})
\to T_{P/J^d(\wt C)}\to,$$
we get $\SF'=T_{P/J^d(\wt C)}(-D_1-D_2)=\SO_P$ and hence
$$h^0(S^2T_X(log))\le \frac{g(g-1)}{2}+h^0(T_{J^d(\wt C)})+1\,=\frac{g(g+1)}{2}.$$
\enddemo

\proclaim{Proposition 3.2} Let $\Cal C/S$ be a flat family of proper curves 
satisfing the assumpations of Proposition 1.3 and such that
$\Cal C_s$ ($s\in\Delta\subset S$) are irreducible curves of one node.    
Let $f:J(\Cal C)\to S$ be the associated family of moduli spaces of torsion
free sheaves of rank $1$ and degree $0$, and $\SL$ be the relative theta line 
bundle on $J(\Cal C)/S$. Then for any integer $k>0$ and $s\in \Delta$
\roster\item $H^0(T_{J(\Cal C_s)}(log))@>\cup c_1(\SL^k_s)>>H^1(\SO_{J(\Cal C_s)})$ 
is an isomorphism.
\item $\mu_{\SO_{J(\Cal C_s)}}=0$ and 
$H^0(\SD^2_{J(\Cal C_s)^{\dagger}}(\SL_s^k))\cong\Bbb C.$\endroster\endproclaim

\demo{Proof} From the discussions in Section 1, the log structure on 
$J(\Cal C_s)$ induced by $logf^{-1}(\Delta)$ is a logarithmic structure in the
sense of [KN], thus we can use our Proposition 3.1 and Lemma 3.4, the (1) is
a corollary of Proposition 3.1.

The claim $\mu_{\SO_{J(\Cal C_s)}}=0$ is equivalent to that 
$$h^0(\SD^2_{J(\Cal C_s)^{\dagger}}(\SO_{J(\Cal C_s)}))=
h^0(\SD^1_{J(\Cal C_s)^{\dagger}}(\SO_{J(\Cal C_s)}))+h^0(S^2T_{J(\Cal C_s)}(log)),$$
which is true for $s\in S\ssm\Delta$ (see [We]). Therefore, by using the semicontinuity
and Lemma 3.4, it is true for all $s\in S$ if we remark that 
$h^0(\SD^1_{J(\Cal C_s)^{\dagger}}(\SO_{J(\Cal C_s)}))$ is constant for all 
$s\in S$ since the
canonical exact sequence 
$$0\to\SO_{J(\Cal C_s)}\to\SD^1_{J(\Cal C_s)^{\dagger}}(\SO_{J(\Cal C_s)})\to
T_{J(\Cal C_s)}(log)\to 0$$
is splitting by Proposition 2.2 and $c_1(\SO_{J(\Cal C_s)})=0$. By using again
Proposition 2.2 and the above (1), we know that $\mu_{\SL^k_s}=\,-\cup c_1(\SL_s^k)$ 
is injective. Hence 
$$H^0(\SD^2_{J(\Cal C_s)^{\dagger}}(\SL_s^k))=
H^0(\SD^1_{J(\Cal C_s)^{\dagger}}(\SL_s^k))\cong\Bbb C.$$\enddemo

\proclaim{Theorem 3.1} Let $f:J(\Cal C)\to S$ be the family of generalized
Jacobians in Proposition 3.2, $Y=f^{-1}(\Delta)$ and $\SL$ be the relative 
theta line bundles. Then there exists a symbol
$$\rho:T_S(log\Delta)\to f_*S^2T_{J(\Cal C)/S}(logY)$$
such that the following two conditions hold\roster
\item $\mu_{\SL^k}\cdot\rho+\kappa_{J(\Cal C)/S}=0,$
\item $f_*T_{J(\Cal C)/S}(logY)@>\cup c_1(\SL^k)>>R^1f_*\SO_{J(\Cal C)}$ is 
an isomorphism.\endroster
In particular, there exists a unique projective logarithmic heat operator
$$\wt H:T_S(log\Delta)\to f_*\SW_{J(\Cal C)/S}(\SL^k)/\SO_S$$
such that $\rho_{\wt H}=\rho$, and thus there exists a projective
logarithmic connection on $f_*\SL^k$.\endproclaim

\demo{Proof} It is clear that we only need to check (1) since (2) has been shown
in Proposition 3.2, namely, we need to find a solution of 
$\mu_{\SL^k}\cdot\rho+\kappa_{J(\Cal C)/S}=0.$
By (2), we have the isomorphism 
$$f_*T_{J(\Cal C)/S}(logY)\otimes f_*T_{J(\Cal C)/S}(logY)
@>\cup c_1(\SL^k)>> R^1f_*T_{J(\Cal C)/S}(logY).$$
Let $\rho=(\cup c_1(\SL^k))^{-1}\circ\kappa_{J(\Cal C)/S}:
T_S(log\Delta)\to f_*T_{J(\Cal C)/S}(logY)\otimes f_*T_{J(\Cal C)/S}(logY)$,
which, over the open set $S\ssm\Delta$, is a map into $f_*S^2T_{J(\Cal C)/S}(logY)$ 
(see \S2.3.8 of [GJ] or [We]), thus it is a map  
$$\rho=(\cup c_1(\SL^k))^{-1}\circ\kappa_{J(\Cal C)/S}:
T_S(log\Delta)\to f_*S^2T_{J(\Cal C)/S}(logY).$$
By Proposition 2.2, $\mu_{\SL^k)}=\,-\cup c_1(\SL^k)$ and $\rho$ is a solution 
of $\mu_{\SL^k}\cdot\rho+\kappa_{J(\Cal C)/S}=0.$ 
\enddemo

\enddocument